\documentclass[11pt]{article}
\usepackage[utf8]{inputenc}
\usepackage[english]{babel}
\usepackage[T1]{fontenc}
\usepackage{float}
\usepackage{listings}
\usepackage{amsmath}
\usepackage{amsfonts}
\usepackage{amssymb}
\usepackage{amsthm}
\usepackage{enumerate}
\usepackage{anysize}
\usepackage{graphicx}
\usepackage{xcolor}
\usepackage{color}
\usepackage{subcaption}
\usepackage{fancyhdr}
\usepackage{makeidx}
\usepackage{hyperref}
\usepackage{appendix}
\usepackage{tikz}
\usepackage[matrix,arrow]{xy}
\usetikzlibrary{arrows} 
\usepackage{multirow}
\usepackage{url}
\usepackage[square,numbers]{natbib}
\usepackage{subcaption}
\usepackage{stackrel}
\marginsize{2cm}{2cm}{0cm}{2cm}
\usepackage{lipsum}

\definecolor
{gray75}{gray}{0.75}
\definecolor
{gray85}{gray}{0.85}
\numberwithin{equation}{section}
\newtheorem{teo}{Theorem}[section]

\newtheorem{obs}{Remark}[section]

\begin{document}

\title{\textbf{Determining parameters giving different growths of a new Glioblastoma differential model}}
\author{ A. Fernández-Romero\footnote{\small{Dpto. Ecuaciones Diferenciales y Análisis Numérico, Facultad de Matemáticas and IMUS, Universidad de Sevilla. Sevilla, Spain; e-mail: \texttt{afernandez61@us.es} \& \texttt{guillen@us.es} \& Corresponding author: \texttt{suarez@us.es}.} }$\;$, F. Guillén-González$^{*}$ \footnote{ORCID: 0000-0001-5539-5888}  $\;$and A. Suárez$^{*}$ \footnote{ORCID: 0000-0002-6407-7758}.
%\blfootnote{\textbf{Keywords:} Glioblastoma, numerical simulation, mathematical oncology, parameter estimation.}
}
%\author{ A. Fernández-Romero$^{2}$, F. Guillén-González$^{2}$\footnote{ORCID: 0000-0001-5539-5888}, A. Suárez$^{1\;2}$\footnote{ORCID: 0000-0002-6407-7758}.\\ \small{$^{1}$Corresponding author.}\\
%	\small{$^{2}$Dpto. Ecuaciones Diferenciales y Análisis Numérico,}\\
%	\small{Facultad de Matemáticas, Universidad de Sevilla. Sevilla, Spain.}\\
%	\small{\texttt{afernandez61@us.es, guillen@us.es, suarez@us.es}}	
%}
\date{}
\maketitle
\section*{\centering{\textbf{Abstract}}}
In this paper we analyse a differential system related to a Glioblastoma growth. Using numerical simulations, we prove that model captures different kind of growth changing adequately the parameters of the model. Firstly, we make an adimensional study in order to reduce the number of parameters. Later, we detect the main parameters determining either different width of the ring formed by proliferative cells around necrotic ones or different regular/irregular behaviour of the tumor surface.
\\
\\
\textbf{Mathematics Subject Classification.} $35\text{M}10,\;35\text{Q}92,\;92\text{B}05,\;92\text{C}17,\;92\text{D}25$.\\
\textbf{Keywords:} Glioblastoma, numerical simulation, mathematical oncology, parameter estimation.\\
\textbf{Funding}: The authors were supported by PGC2018-098308-B-I00 (MCI/AEI/FEDER, UE).
%\tableofcontents
\section{Introduction}
Glioblastoma (GBM) is one of the more lethal brain tumors with a survival from $6$ to $14$ months for those patients who receive an standard care \cite{Davis_2016,Ostrom_2014}. Due to this fact, GBM is studied with a high interest from the oncology community even from a mathematical point of view (see \cite{Alfonso_2017,Baldock_2013,Protopapa_2018} and references therein). Some works such as \cite{Rockne_2009} are dedicated to model the GBM growth and control its postoperative treatment whereas some PDEs models have been used in order to analyse the proliferation and cellular invasion inside the brain, see \cite{Swanson_2008, Swanson_2000, Tracqui_1995} and the references cited therein.
\\

\begin{figure}[H]
	\includegraphics[width=6cm, height=6cm]{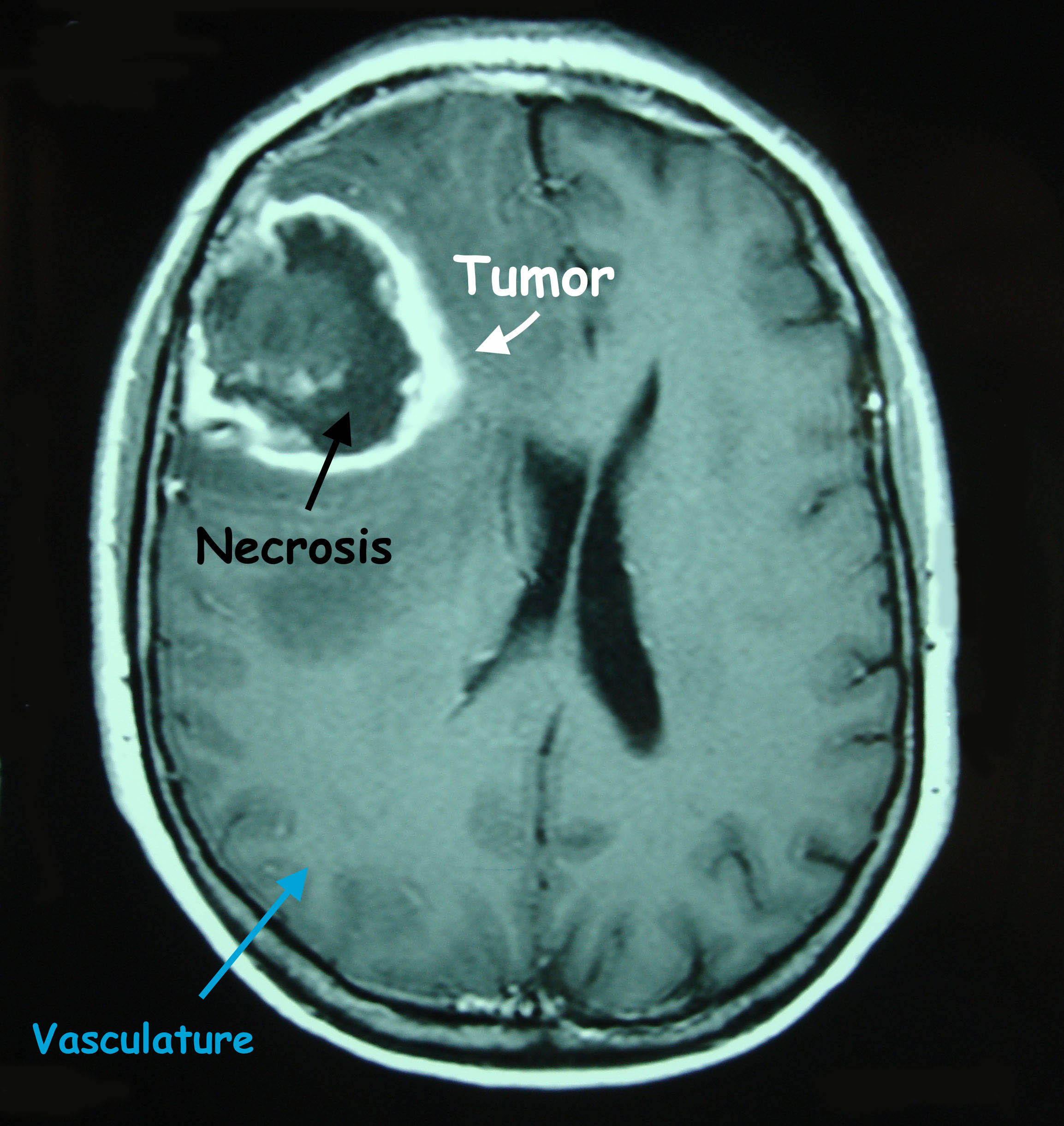}
	\centering
	\caption{MRI of a GBM areas showing contrast enhancement\protect\footnotemark~.}
	\label{MRI}
\end{figure}

The magnetic resonance images (MRIs) of some GBMs show a necrotic area in the centre surrounded by a white ring as we see in Figure $\ref{MRI}$. This ring corresponds to the zone with proliferative tumor and it is also an indicator of areas with low vasculature.
\\

\footnotetext{%Tratamiento y resultados del Glioblastoma multiforme en nuestro centro. Periodo 2010-2014. Máster en avances en radiología diagnóstica y terapéutica en medicina física. Universidad de Granada. Sara Rodríguez Pavón. 
	\url{https://pdfs.semanticscholar.org/7d7b/2f5f038cf961be42c789db6a8dffa8637733.pdf}}

\begin{figure}[H]
	\includegraphics[width=10cm, height=8cm]{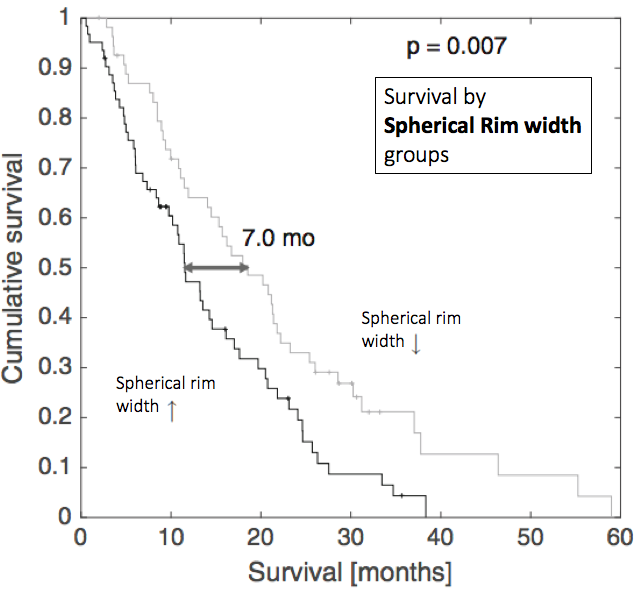}
	\centering
	\caption{Survival vs the spherical rim width of GBM, \cite{Julian_2016}.}
	\label{fig:fig1}
\end{figure}

The necrotic area is an important characteristic of GBMs since it can determine the volume of the GBM and its prognosis in relation to mortality. According to that, an experimental study relating the ring width of proliferative tumor around the necrosis and its mortality is shown in \cite{Julian_2016} (see Figure $\ref{fig:fig1}$).

%\begin{figure}[H]
%	\includegraphics[width=0.6\linewidth]{output2}
%	\centering
%	\caption{Survival vs the ring width-volume of GBM.}
%	\label{fig:fig1}
%\end{figure}

The study of \cite{Julian_2016} concludes that tumors with slim ring have better prognostic, specifically $7$ months of more survival, than tumors with thick ring. Other way to understand this study is based in the amount of necrosis, since tumors with slim ring have more amount of necrosis than those with thick ring.
\\

\begin{figure}[H]
	\includegraphics[width=10cm, height=8cm]{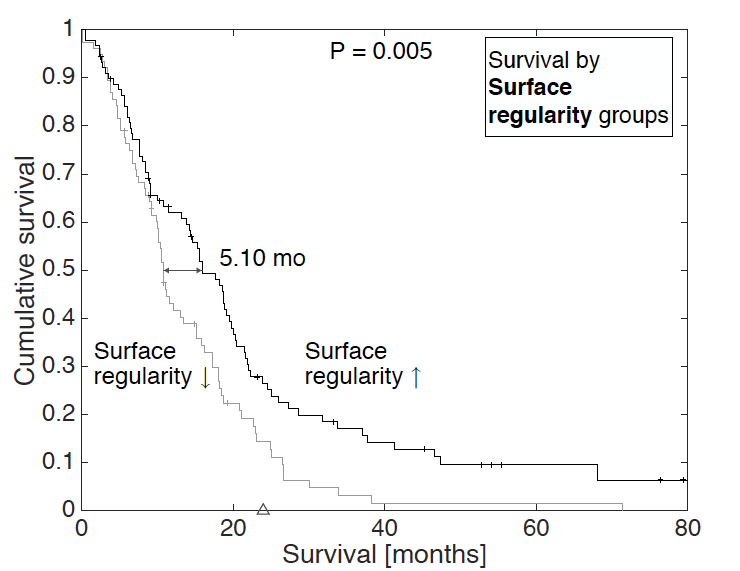}
	\centering
	\caption{Survival vs the regularity surface of GBM, \cite{Victor_2018}.}
	\label{fig:fig4}
\end{figure}

Another relevant aspect of a GBM observed in the MRIs is the regularity of the tumor surface. In \cite{Victor_2018}, the authors made an experimental study about the survival of patients in relation to the regular or irregular surface growth, of the GBM. Indeed, Figure $\ref{fig:fig4}$ shows that tumors with a regular surface have better prognostic, more than $5$ moths of survival, than tumor with irregular surface. 

Once we have presented these two main characteristics of the GBM, our goal in this paper is to use a differential model including the variables for interactions between proliferative tumor, necrosis and vasculature, in order to capture the two phenomena presented above through numerical simulations and to detect which parameters are more important in each kind of GBM growth. Thus, we will study two different growths: the first one consists of computing the so-called tumor-ring via ratio proliferative tumor/necrosis and the second one is about to detect regular vs irregular growth of the tumor depending on vasculature. For these studies, we have designed two coefficients depending on: the density for the ring widht, defined by "ring quotient (RQ), and the area for the regularity surface, defined by "surface quotient (SQ). Using these coefficients and changing the value of the parameters of the model, we obtain the relevancy of the parameters in the different tumor growths.
\\

The structure of the paper is the following: In Section $\ref{modelo}$, we present the model used and the main mathematical results obtained in the previous work \cite{Romero2_2020}. In Section $\ref{adimensionalizacion}$, an adimesionalization of the model is showed to reduce the number of parameters. Next, in Section $\ref{anillo}$, we make the study of the ring width in relation to the parameters. Section $\ref{regularidad}$ is devoted to study the regularity surface of the GBM with respect to the parameters. Finally, in Section $\ref{conlcusion}$, we discuss and summarize our results.

\section{The Model}\label{modelo} 

Here, we present the nonlinear diffusion model, studied in \cite{Romero2_2020}, that we will use along the paper:
\begin{equation}\label{probOriginal}
\left\{\begin{array}{ccl}
\dfrac{\partial T}{\partial t} -\nabla\cdot\left(\left(\kappa_1\;P\left(\Phi,T\right)+\kappa_0\right)\nabla\;T\right)& = & f_1\left(T,N,\Phi\right)\quad\text{in}\;\;\left(0,T_f\right]\times\Omega\\ %\;\;\text{en}\;\;\left(0,T_f\right)\times\Omega\\
&&\\
\dfrac{\partial N}{\partial t}& = & f_2\left(T,N,\Phi\right)\quad\text{in}\;\;\left(0,T_f\right]\times\Omega\\%\;\;\text{en}\;\;\left(0,T_f\right)\times\Omega \\
&&\\
\dfrac{\partial \Phi}{\partial t} & = &f_3\left(T,N,\Phi\right)\quad\text{in}\;\;\left(0,T_f\right]\times\Omega%\;\;\text{en}\;\;\left(0,T_f\right)\times\Omega \\
\\
%\hspace{0.5cm}\dfrac{\partial T}{\partial \text{n}}\Bigg\vert_{\partial\Omega}=0\;\;\text{en}\;\;\left(0,T_f\right)\times\partial\Omega\\
%\\

\end{array}\right.\end{equation}
\\
with non-flux boundary condition
\begin{equation}\label{condifronte}
%\dfrac{\partial T}{\partial n}\Bigg\vert_{\partial\Omega}=0\;\;\text{on}\;\;\left(0,T_f\right)\times\partial\Omega
\nabla T\cdot n=0\;\;\text{on}\;\;\left(0,T_f\right)\times\partial\Omega
\end{equation} 
where $n$ is the outward unit normal vector to $\partial\Omega$, and initial conditions
\begin{equation}\label{condinicio}
T\left(0,\cdot\right)=T_0(x),\;N\left(0,\cdot\right)=N_0(x),\;\Phi\left(0,\cdot\right)=\Phi_0(x)
\;\;\text{in}\;\;\Omega.
\end{equation}

The domain $\Omega\subset\mathbb{R}^2$ or $\mathbb{R}^3$ is bounded and regular, $T_f>0$ is the final time, and $T(x,t), N(x,t)$ and $\Phi(x,t)$ represent tumor (proliferative) and necrotic densities and the vasculature concentration at the point $x\in\Omega$ and the time $t>0$, respectively. The nonlinear reactions terms of $\left(\ref{probOriginal}\right)$ are defined by

%\begin{equation}\label{funciones}
%\left\{\begin{array}{ccl}
%f_1\left(T,N,\Phi\right) &:=&\rho\;T\;P\left(\Phi,T\right) \left(1-\dfrac{T+N+\Phi}{K}\right)-\alpha\;T\;\sqrt{1-P\left(\Phi,T\right)^2}-\beta_1\; 
%N\;T,\\
%\\
%f_2\left(T,N,\Phi\right) &:=& \alpha\;T\;\sqrt{1-P\left(\Phi,T\right)^2}+\beta_1\; N\;T +\delta\;T\;\Phi +\beta_2\;N\;\Phi,\\
%\\
%f_3\left(T,N,\Phi\right) &:=& \gamma\;T\;\sqrt{1-P\left(\Phi,T\right)^2}\;\dfrac{\Phi}{K}\left(1-\dfrac{T+N+\Phi}{K}\right)-\delta\;T\;\Phi-\beta_2\;N\;\Phi.
%\end{array}\right.
%\end{equation}

\begin{subequations}\label{funciones}
	\begin{align}
	f_1\left(T,N,\Phi\right) &:=\underbrace{\rho\;T\;P\left(\Phi,T\right) \left(1-\dfrac{T+N+\Phi}{K}\right)}_{\text{Tumor growth}}-\underbrace{\alpha\;T\;\sqrt{1-P\left(\Phi,T\right)^2}}_{\text{Hypoxia}}-\underbrace{\beta_1\; 
	N\;T,}_{\substack{\text{Tumor destruction}\\\text{ by necrosis}}}\label{eqT}\\
	\nonumber\\
	f_2\left(T,N,\Phi\right) &:= \underbrace{\alpha\;T\;\sqrt{1-P\left(\Phi,T\right)^2}}_{\text{Hypoxia}}+\underbrace{\beta_1\; 
		N\;T,}_{\substack{\text{Tumor destruction}\\\text{ by necrosis}}} +\underbrace{\delta\;T\;\Phi}_{\substack{\text{Vascular destruction}\\ \text{by tumor}}} +\underbrace{\beta_2\;N\;\Phi.}_{\substack{\text{Vasculature destruction}\\ \text{by necrosis }}}\label{eqN}\\
	\nonumber\\
	\nonumber f_3\left(T,N,\Phi\right) &:= \underbrace{\gamma\;T\;\sqrt{1-P\left(\Phi,T\right)^2}\;\dfrac{\Phi}{K}\left(1-\dfrac{T+N+\Phi}{K}\right)}_{\text{Vasculature growth}}-\underbrace{\delta\;T\;\Phi}_{\substack{\text{Vascular destruction}\\ \text{by tumor}}}\\
	\nonumber\\
	&-\underbrace{\beta_2\;N\;\Phi.}_{\substack{\text{Vasculature destruction}\\ \text{by necrosis }}} \label{eqPhi}
	\end{align}
\end{subequations}

The function $P\left(\Phi,T\right)$ will be a ratio between vasculature and tumor plus vasculature, defined as follows
\begin{equation}\label{funcionP}
P\left(\Phi, T\right)=
\dfrac{\Phi_+}{\left(\dfrac{\Phi_++K}{2}\right)+T_+}
\end{equation}
with $T_+=\max\{0,T\}$ and similar to $\Phi_+$. Notice that the vasculature volume fraction $P\left(\Phi,T\right)$ is a continuous function in $\mathbb{R}^2$, satisfying the pointwise estimates 
%$$0\leq P\left(\Phi,T\right)\leq2\quad\forall\left(T,\Phi\right)\in\mathbb{R}^2,$$
$$0\leq P\left(\Phi,T\right)\leq1\quad\forall\left(T,\Phi\right)\in\left[0,K\right]\times\left[0,K\right]$$
and $P\left(\Phi,T\right)=0$ for $\Phi=0$. On the other hand, the factor $\sqrt{1-P\left(\Phi,T\right)^2}$ acting in the hypoxia term can be seen as a volume fraction measuring the lack of vasculature, and it has  the same pointwise estimates that $P\left(\Phi,T\right)$.
\\

The parameters $\kappa_1$, $\kappa_0$, $\rho$, $\rho$, $\alpha$, $\beta_1$, $\beta_2$, $\gamma$, $\delta$, $K>0$ in $\left(\ref{probOriginal}\right)$ have the following description \cite{Klank_2018,Alicia_2015, Alicia_2012}:

\begin{table}[H]
	\centering
	\begin{tabular}{c|c|c}
		%\hline
		\textbf{Variable} & \textbf{Description} & \textbf{Value} \\
		\hline
		$\kappa_1$	& Anisotropic  speed diffusion & $cm^2/\text{day}$     \\
		\hline
		$\kappa_0$	& Isotropic speed diffusion & $cm^2/\text{day}$     \\
		\hline
		$\rho$	& Tumor proliferation rate  &  $\text{day}^{-1}$    \\
		\hline
		$\alpha$ & Hypoxic death rate by persistent anoxia & $cell/\text{day}$    \\
		\hline
		$\beta_1$	& Change rate from tumor to necrosis &   $\text{day}^{-1}$    \\
		\hline
		$\beta_2$	& Change rate from vasculature to necrosis &   $\text{day}^{-1}$    \\
		\hline
		$\gamma$	& Vasculature proliferation rate  &  $\text{day}^{-1}$    \\
		\hline
		$\delta$	& Vasculature destruction by tumor action  &  $\text{day}^{-1}$    \\
		\hline
		$K$	& Carrying capacity & $\text{cell}/\text{cm}^3$ \\
		%\hline    
	\end{tabular}
	\caption{\label{parametros} Parameters.}
\end{table}

In $\left(\ref{probOriginal}\right)$, we model a movement of tumor by diffusion with higher velocity in zones with more vasculature. Thus, this effect is expressed by function $P\left(\Phi,T\right)$ given in $\left(\ref{funcionP}\right)$, measuring the quotient between the amount of vasculature and the amount of vasculature and tumor together. Moreover, in $\left(\ref{eqT}\right)$ the velocity of the tumor growth is proportional to $P(\Phi,T)$ since vasculature supplies nutrients and oxygenation to cells. Conversely, tumor also decreases with the lack of vasculature in the hypoxia death rate of tumor cells. Furthermore, a destruction of tumor by necrosis is considered by the term $-\beta_1\;N\;T$. In $\left(\ref{eqPhi}\right)$, we have a logistic growth term since vasculature needs space to growth. However, the speed of the vasculature growth depends on tumor through the term $T\;\sqrt{1-P\left(\Phi,T\right)^2}$ showing that there will not be growth of vasculature in absence of tumor and this speed increases when vasculature decreases. Furthermore, a destruction of vasculature by tumor and necrosis are considered by the terms $-\delta\;\Phi\;T$ and $-\beta_2\;N\;\Phi$, respectively. 
\\

Finally in equation $\left(\ref{eqN}\right)$ we observe as necrosis is formed by the sum of death terms of tumor and vasculature.
\\

Let us point out that, recently, in \cite{Victor_2020}, the authors have proposed a mathematical model, simpler than $\left(\ref{probOriginal}\right)$, %the one presented above, 
for a GBM growth to quantify only the tumor ring and the a relation with the survival showed in Figure $\ref{fig:fig1}$. We have completed the model of \cite{Victor_2020} in order to not only capture the ring width but also the regularity surface of the GBM. 
\\

Another advantage of $\left(\ref{probOriginal}\right)$ over the model of \cite{Victor_2020} is the presence of the vasculature as an additional variable which is essential in the study of the regularity surface, as we will see in Section $\ref{regularidad}$, since the amount and spatial distribution of vasculature determine the growth of tumor. Moreover, the introduction of vasculature would allow the application of chemical therapies in the model because this type of therapy is driven by the vasculature.
\\

In a previous work, see \cite{Romero2_2020}, problem $\left(\ref{probOriginal}\right)$-$\left(\ref{condinicio}\right)$ has been studied mathematically (from analysis to numerics), where the following results were obtained:

\begin{teo}[Existence of global in time Solution of Problem $\left(\ref{probOriginal}\right)$-$\left(\ref{condinicio}\right)$]\label{solucion_Problema}
	Given $T_0\in L^\infty\left(\Omega\right)$ and $N_0,\Phi_0\in H^1\left(\Omega\right)\cap L^\infty(\Omega)$ satisfying $0\leq T_0(x),N_0(x),\Phi_0(x)\leq K$, a.e.$\;x\in\Omega$, then problem $\left(\ref{probOriginal}\right)$-$\left(\ref{condinicio}\right)$ has a solution $\left(T, N,\Phi\right)$ such that 
	$0\leq T,\Phi\leq K$, $0\leq N\leq C\left(T_f\right)$
	with $C\left(T_f\right)$ a positive constant which depends exponentially on the final time $T_f>0$ and the carrying capacity $K$. Indeed,
	$T\in L^\infty\left(0,T_f;L^2\left(\Omega\right)\right)\cap L^2\left(0,T_f;H^1\left(\Omega\right)\right)$, $T_t\in L^2\left(0,T_f;\left(H^1\left(\Omega\right)\right)'\right)$, 
	$N,\Phi\in{L}^\infty\left(0,T_f;H^1\left(\Omega\right)\right)$, $N_t,\Phi_t\in{L}^2\left(0,T_f;L^2\left(\Omega\right)\right)$
	and they satisfy
	%satisfies the variational formulation
	$$\displaystyle \int_{0}^{T_f} \langle T_t,v\rangle_{\left(H^1\left(\Omega\right)\right)'}\;dt+\int_{0}^{T_f}\int_\Omega \left(\kappa_1\;P\left(\Phi,T\right)+\kappa_0\right)\nabla T\cdot\nabla v\;dx\;dt=\int_{0}^{T_f}\int_\Omega f_1\left(T,N,\Phi\right)\;v\;dx\;dt,$$
	$\forall v\in L^2\left(0,T_f;H^1\left(\Omega\right)\right)$ and %\quad\text{where the space}\;\;W_2\;\;\text{is defined by } \ref{W}$$
	%	and
	%	$N,\Phi\in{L}^\infty\left(0,T_f;H^1\left(\Omega\right)\right),\;\;N_t,\Phi_t\in{L}^2\left(0,T_f;L^2\left(\Omega\right)\right)$
	%	verifying
	
	$$\left\{\begin{array}{rl}
	\displaystyle N_t=f_2\left(T,N,\Phi\right)&\\
	&\quad\text{a.e. in}\;\; %\left(t,x\right)\in
	\left(0,T_f\right)\times\Omega\\
	\displaystyle \Phi_t=f_3\left(T,N,\Phi\right)&
	\end{array}\right.$$
	%$$\displaystyle N_t=f_2\left(T,N,\Phi\right)$$
	%$$\displaystyle \Phi_t=f_3\left(T,N,\Phi\right)$$
	%a.e. $\left(t,x\right)\in\left(0,T_f\right)\times\Omega$ 
	and the boundary and initial conditions $\left(\ref{condifronte}\right)$ and $\left(\ref{condinicio}\right)$ are satisfied by $T$ and $(T,N,\Phi)$, respectively.
\end{teo}

\begin{teo}[Long time behaviour]\label{lemaFi_0}
	Given $\epsilon>0$ and a solution $\left(T,N,\Phi\right)$ of $\left(\ref{probOriginal}\right)$-$\left(\ref{condinicio}\right)$, if there exists $\widetilde{\Omega}\subset\Omega$ with $\mid\widetilde{\Omega}\mid>0$ such that $0<\epsilon\leq N_0\left(x\right)$ a.e.$\;x\in \widetilde \Omega$, one has $\Phi\left(t,x\right)\rightarrow0$ when $t\rightarrow+\infty$ a.e. $x\in \widetilde \Omega$. In addition, if $\delta\geq\dfrac{\gamma}{K}$, it also holds that $T\left(t,x\right)\rightarrow0$ when $t\rightarrow+\infty$ a.e.$\;x\in \widetilde \Omega$. Moreover, there exists $N_{\max}\ge\|N_0\|_{L^\infty\left(\Omega\right)}$ such that
	$N\left(t,x\right)\leq N_{\max}$ a.e. $\left(t,x\right)\in\left(0,+\infty\right)\times\Omega$.
	
\end{teo}

In addition, the numerical scheme that we will use in this work, have been presented in \cite{Romero2_2020},  proving that this scheme preserves the pointwise and energy estimates showed in Theorem $\ref{solucion_Problema}$.

\section{Adimensionalization}\label{adimensionalizacion}

Before showing the numerical simulations related to the different GBM growths, we make a study about the parameters, simplifying and presenting only the simulations according to the relevant adimensional parameters.
\\

The first study will depend on the carrying capacity, the parameter $K>0$. We consider the change of variables $\widetilde{T}=\dfrac{T}{K}$, $\widetilde{N}=\dfrac{N}{K}$ and $\widetilde{\Phi}=\dfrac{\Phi}{K}$ passing the normalized capacity to $1$.
\\

For the second adimensionalization, we consider the parameters $\kappa_0$ and $\rho$ as our point of study. Since $\rho$ corresponds to tumor proliferation rate, it is related with the time while the diffusion parameter $\kappa_0$ is related to the spatial variable. Thus, we can make the following change of the independent variables:
\begin{equation}\label{cambios_variable}
\left\{
\begin{array}{lcl}
s=\rho\;t&\Rightarrow& ds=\rho\;dt,\\
\\
y=\sqrt{\dfrac{\rho}{\kappa_0}}\;x&\Rightarrow& dy=\sqrt{\dfrac{\rho}{\kappa_0}}\;dx.
\end{array}
\right.
\end{equation}

Applying these changes, our system $\left(\ref{probOriginal}\right)$ becomes to 

\begin{equation}\label{prob_K_rho}
\left\{\begin{array}{ccl}
\dfrac{\partial \widetilde{T}}{\partial s} -\nabla\cdot\left(\left(\dfrac{\kappa_1}{\kappa_0}\;P\left(\widetilde{\Phi},\widetilde{T}\right)+1\right)\nabla\;\widetilde{T}\right)& = & \widetilde{f_1}\left(\widetilde{T},\widetilde{N},\widetilde{\Phi}\right)\\ 
&&\\
\dfrac{\partial \widetilde{N}}{\partial s}& = & \widetilde{f_2}\left(\widetilde{T},\widetilde{N},\widetilde{\Phi}\right)\\
&&\\
\dfrac{\partial \widetilde{\Phi}}{\partial s} & = &\widetilde{f_3}\left(\widetilde{T},\widetilde{N},\widetilde{\Phi}\right)\\
\end{array}\right.\end{equation}
where
\begin{equation}\label{funciones_K_rho}
\left\{\begin{array}{lll}
\widetilde{f_1}\left(\widetilde{T},\widetilde{N},\widetilde{\Phi}\right) &=&\widetilde{T}\;P\left(\widetilde{\Phi},\widetilde{T}\right) \left(1-\left(\widetilde{T}+\widetilde{N}+\widetilde{\Phi}\right)\right)-\dfrac{\alpha}{\rho}\;\widetilde{T}\;\sqrt{1-P\left(\widetilde{\Phi},\widetilde{T}\right)^2}-K\;\dfrac{\beta_1}{\rho}\; 
\widetilde{N}\;\widetilde{T},\\
\\
\widetilde{f_2}\left(\widetilde{T},\widetilde{N},\widetilde{\Phi}\right) &=& \dfrac{\alpha}{\rho}\;\widetilde{T}\;\sqrt{1-P\left(\widetilde{\Phi},\widetilde{T}\right)^2}+K\;\dfrac{\beta_1}{\rho}\; \widetilde{N}\;\widetilde{T} +K\;\dfrac{\delta}{\rho}\;\widetilde{T}\;\widetilde{\Phi} +K\;\dfrac{\beta_2}{\rho}\;\widetilde{N}\;\widetilde{\Phi},\\
\\
\widetilde{f_3}\left(\widetilde{T},\widetilde{N},\widetilde{\Phi}\right) &=& \dfrac{\gamma}{\rho}\;\widetilde{T}\;\sqrt{1-P\left(\widetilde{\Phi},\widetilde{T}\right)^2}\;\widetilde{\Phi}\left(1-\left(\widetilde{T}+\widetilde{N}+\widetilde{\Phi}\right)\right)-K\;\dfrac{\delta}{\rho}\;\widetilde{T}\;\widetilde{\Phi}-K\;\dfrac{\beta_2}{\rho}\;\widetilde{N}\;\widetilde{\Phi}.
\end{array}\right.
\end{equation}

Hence, we can rewrite the rest of the dimensionless parameters as follows:

\begin{table}[H]
	\centering
	\begin{tabular}{c|c|c|c|c|c|c}
		%\hline
		\textbf{Dimensionless parameter} & $\kappa_1^*$&$\alpha^*$  &$\beta_1^*$	 &	$\beta_2^*$	 &	$\gamma^*$&	$\delta^*$ \\
		\hline
		 \textbf{Original parameter} &$\vspace{2cm}\dfrac{\kappa_1}{\kappa_0}$   & $\dfrac{\alpha}{\rho}$  & $K\;\dfrac{\beta_1}{\rho}$  & $K\;\dfrac{\beta_2}{\rho}$ & $\dfrac{\gamma}{\rho}$ & $K\;\dfrac{\delta}{\rho}$  
		 \vspace{-2cm}
	\end{tabular}
	\caption{\label{parametros_K_rho} Dimensionless parameters.}
\end{table}

Thus, we have reduced our model in three parameters: $\kappa_0$, $\rho$ and $K$. Moreover, with this simplification, we could obtain the same conclusions depending on every parameter without the necessity to simulate the growth for different $\kappa_0$ and/or $\rho$ since the increase or decrease of $\kappa_0$ and/or $\rho$ will be understood as a increase or decrease of the other parameters.

\begin{obs}
	To simplify the notation, we drop $'*'$ and $'\;\textasciitilde\;'$ along the paper: $s=t$, $y=x$, $\kappa_1^*=\kappa_1$, $\alpha^*=\alpha$, $\beta_i^*=\beta_i$ for $i=1,2$, $\gamma^*=\gamma$, $\delta^*=\delta$, $\widetilde{T}=T$, $\widetilde{N}=N$, $\widetilde{\Phi}=\Phi$ and $\widetilde{f}_i=f_i$ for $i=1,2,3$.
\end{obs}

To get the numerical simulations we will work with an uncoupled and linear fully discrete scheme of $\left(\ref{probOriginal}\right)$-$\left(\ref{condinicio}\right)$ defined in \cite{Romero2_2020} by means of an Implicit-Explicit (IMEX) Finite Difference in time approximation and $P_1$ continuous finite element
with "mass-lumping" in space. The computational domain is fixed to $\Omega=\left(-9,9\right)\times\left(-9,9\right)$ and the final time $T_f=500$. Moreover, this scheme will preserve the pointwise and energy estimates of Theorem $\ref{solucion_Problema}$. In the numerical setting, we construct a structured triangulation $\left\{\mathcal{T}_h\right\}_{h>0}$ of $\overline{\Omega}$ such that $\overline{\Omega}=\displaystyle\bigcup_{\mathcal{K}\in\mathcal{T}_h}\mathcal{K}$ with the partitioning the edges of the boundary of $\Omega$ into $45$ subintervals, corresponding with the mesh size $h=0.4$. Finally, the time step size is chosen as $dt=10^{-3}$.
\\

\begin{figure}[H]
	\includegraphics[width=8cm, height=5cm]{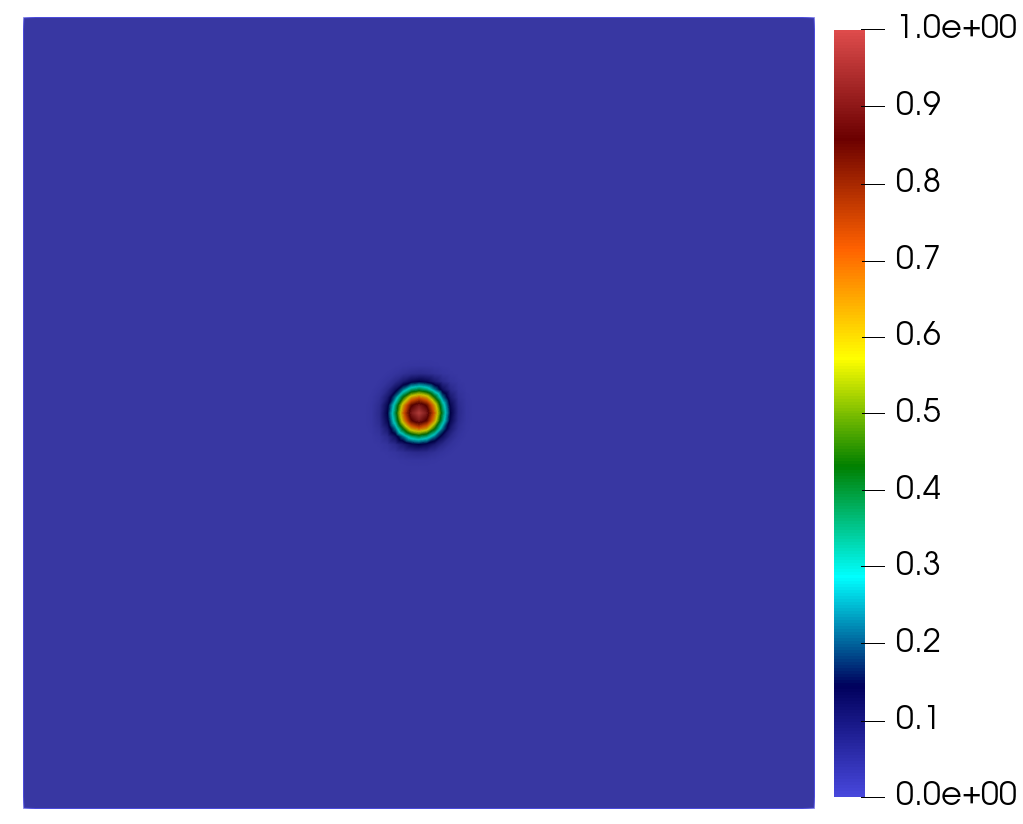}
	\centering
	\caption{Initial tumor.}
	\label{tumor_incial}
\end{figure}

In all the simulations we consider necrosis zero initially and the initial tumor is given by Figure $\ref{tumor_incial}$: 
%The initial conditions for tumor and necrosis in all the simulations will be always the same while the initial condition for vasculature will change depending on the kind of tumor growth studied. Thus, we consider necrosis zero initially and initial tumor given by:

The initial condition for vasculature will change depending on the kind of tumor growth studied

\section{Ring width}\label{anillo}
In order to detect in our model which parameter/s could have more influence in the ring width, we will present simulations according to quantify the tumor-ring with respect to the density of tumor and necrosis. For every simulation, we will change the value of one parameter and checking how the tumor growth changes. 
\\

Since we keep in mind tumor and necrosis, we move the parameters appearing in tumor and necrosis equations, these are, $\kappa_1$, $\alpha$ and $\beta_1$. In all these simulations the value of $\gamma$, $\delta$ and $\beta_2$ are fixed (see Table $\ref{parametros_fijos_razon_tumor_necrosis}$).

\begin{table}[H]
	\centering
	\begin{tabular}{c|c|c|c}
		%\hline
		\textbf{Variable} &$\gamma$ & 	$\delta$ & $\beta_2$ \\
		\hline
		 \textbf{Value} &   $0.255$  &$2.55$    & $2.55$ \\
	\end{tabular}
	\caption{\label{parametros_fijos_razon_tumor_necrosis} Fixed value parameters.}
\end{table}

For the variable parameters, $\kappa_1$, $\alpha$ and $\beta_1$, we will take the following values (see Table $\ref{parametros_variables_razon_tumor_necrosis}$).

\begin{table}[H]
	\centering
	\begin{tabular}{c|c|c|c}
		%\hline
		\textbf{Variable (Fixed value)} & $\kappa_1\;\;\left(55\right)$   &	$\alpha\;\;\left(45\right)$  & 	$\beta_1\;\;\left(27.5\right)$ \\
		\hline
		\textbf{Ranges} & $\left[10,\;100\right]$ & $\left[10,\;100\right]$ &$\left[5,\;50\right]$  \\ 
	\end{tabular}
	\caption{\label{parametros_variables_razon_tumor_necrosis} Variable value parameters.}
\end{table}

Moreover, along this section, we take the initial vasculature defined uniformly in space.
%\\
%
%Finally, the chosen criterion in order to capture the prognosis of tumors will be through the total density of the tumor. That is, we will consider the tumor with more amount of density as the tumor with the worst prognosis.

\subsection{Tumor Ring quotient}

To start with, we show the graphs according to the ratio between proliferative tumor density, $\displaystyle\int_{\Omega}T\;dx$ and total tumor density, $\displaystyle\int_{\Omega}\left(T+N\right)\;dx$, for the different values of $\kappa_1$, $\alpha$ and $\beta_1$ taken in Table $\ref{parametros_variables_razon_tumor_necrosis}$. For that, we define the following "ring quotient" (RQ) coefficient:

\begin{equation}\label{RQ}
0\leq\text{RQ}=\dfrac{\displaystyle\int_{\Omega}T\;dx}{\displaystyle\int_{\Omega}\left(T+N\right)\;dx}\leq1.
\end{equation}

Thus, we can conclude that if RQ is near to zero, tumor ring will be slim due to there exists a high density of necrosis whereas if RQ is close to one, tumor ring will be thick.

%\begin{figure}[H]
%	\includegraphics[width=17cm, height=6cm]{RAZON_TUMOR/EPS/RQ_kappa1.eps}
%	\centering
%	\caption{RQ versus time for $\kappa_1$.}
%	\label{Ring_dif_kappa1}
%\end{figure}
%
%\begin{figure}[H]
%	\includegraphics[width=17cm, height=6cm]{RAZON_TUMOR/EPS/RQ_alpha.eps}
%	\centering
%	\caption{RQ versus time for $\alpha$.}
%	\label{Ring_dif_alpha}
%\end{figure}
%\begin{figure}[H]
%	\includegraphics[width=17cm, height=6cm]{RAZON_TUMOR/EPS/RQ_beta1.eps}
%	\centering
%	\caption{RQ versus time for $\beta_1$.}
%	\label{Ring_dif_beta}
%\end{figure}

\begin{figure}[H]
	\centering
	\begin{subfigure}[b]{0.45\linewidth}
		\includegraphics[width=9cm, height=5cm]{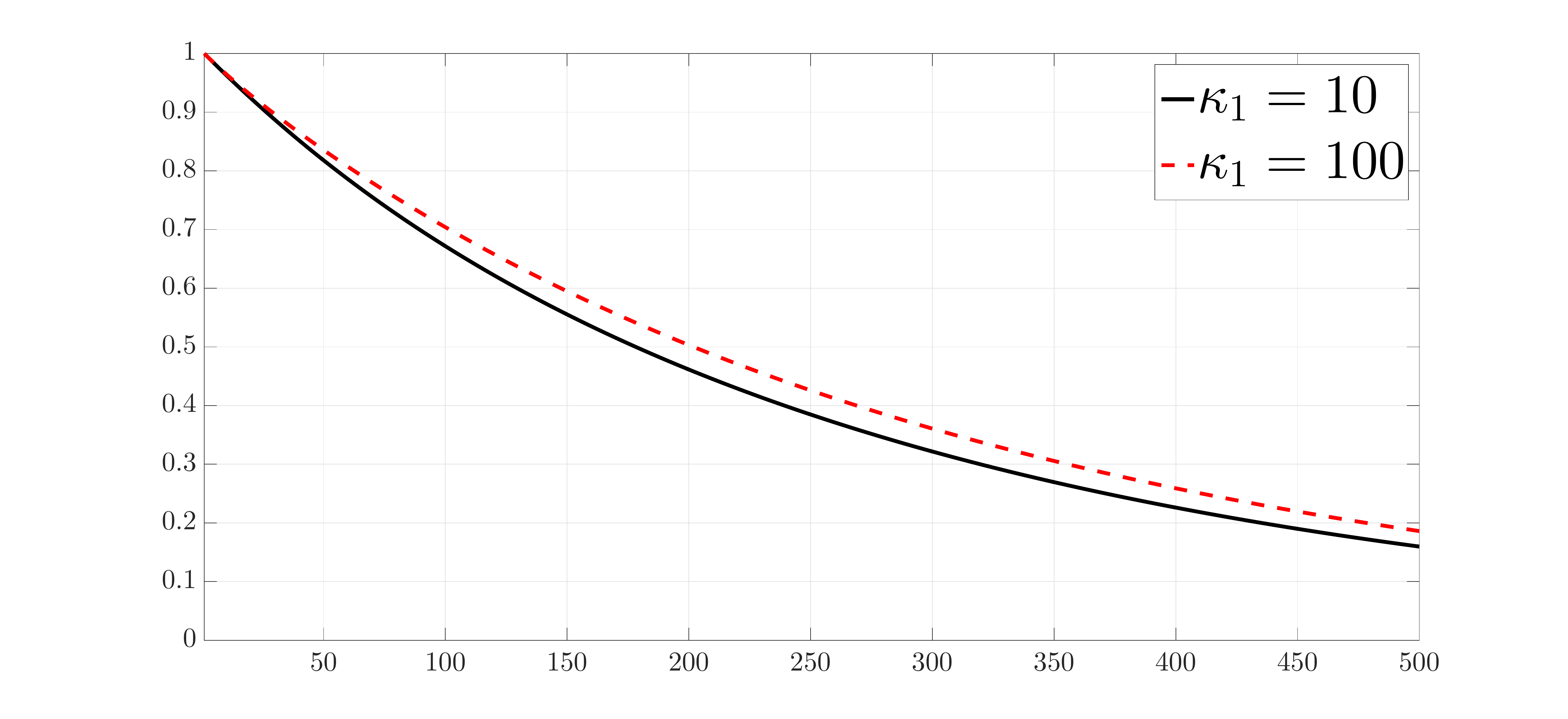}
		\centering
	\caption{RQ versus time for $\kappa_1$.}
	\label{Ring_dif_kappa1}
	\end{subfigure}
	\hspace{1cm}
	\begin{subfigure}[b]{0.45\linewidth}
			\includegraphics[width=9cm, height=5cm]{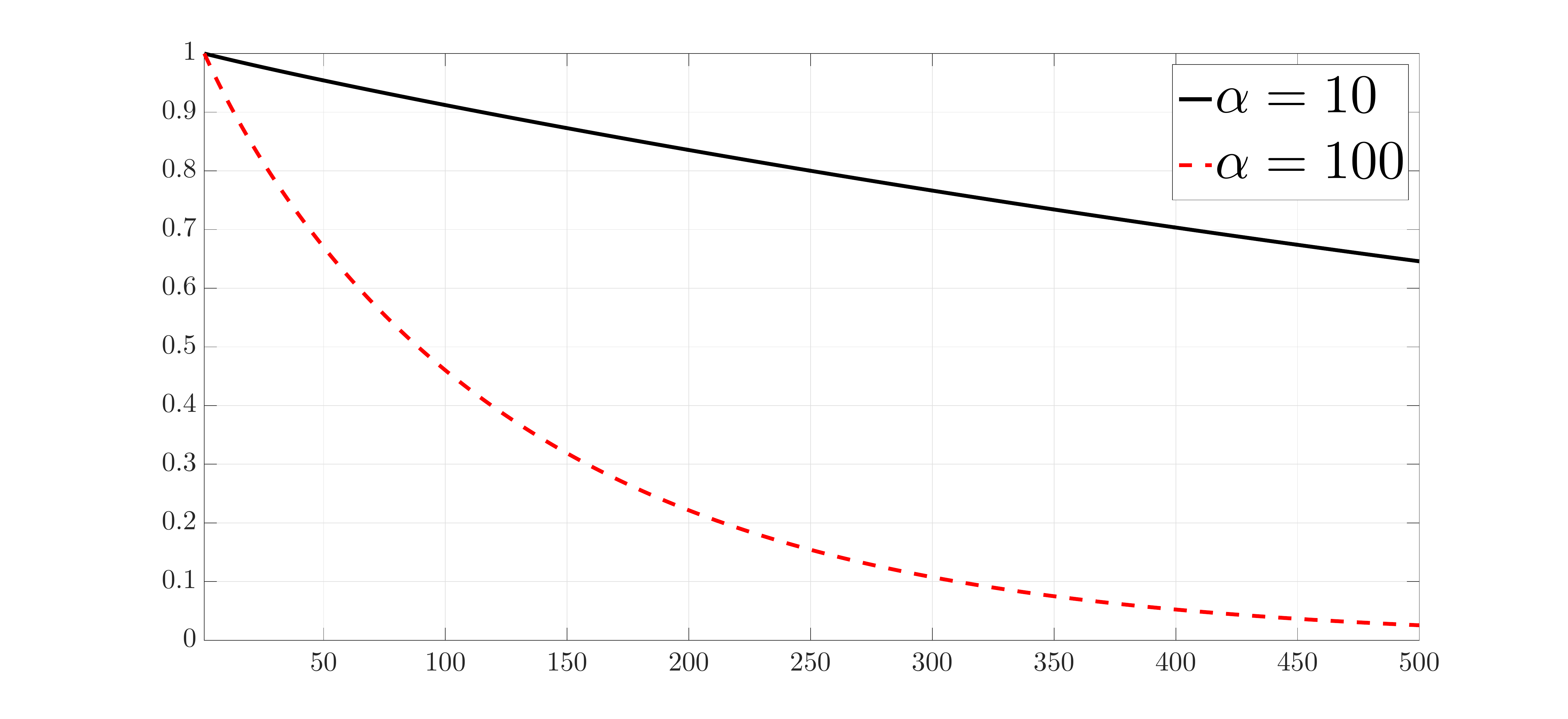}
		\centering
		\caption{RQ versus time for $\alpha$.}
		\label{Ring_dif_alpha}
	\end{subfigure}
\end{figure}
%\vspace{15cm}
\begin{figure}[H]\ContinuedFloat
	\begin{subfigure}[b]{0.45\linewidth}
				\includegraphics[width=9cm, height=5cm]{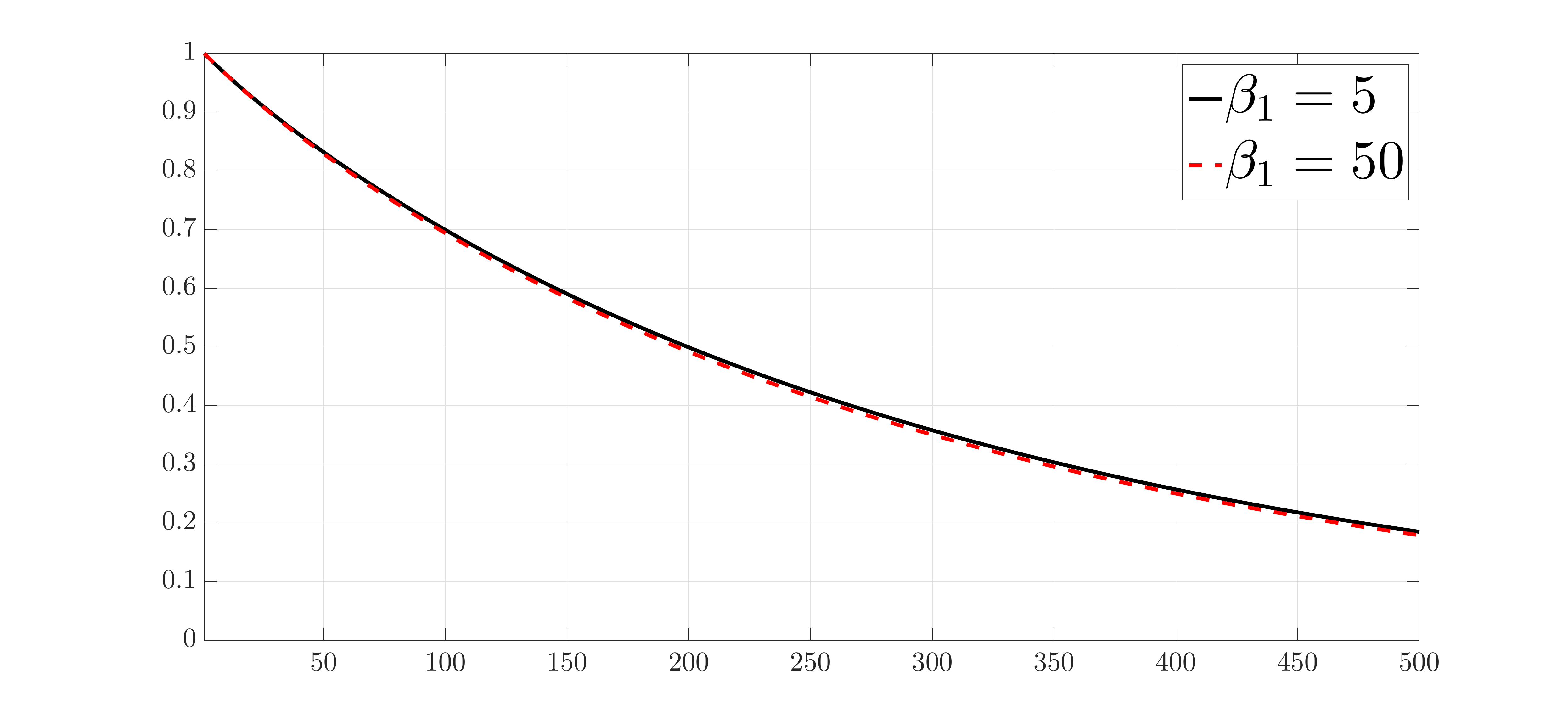}
		\centering
		\caption{RQ versus time for $\beta_1$.}
		\label{Ring_dif_beta}
	\end{subfigure}
	\centering
	\caption{RQ versus time for $\kappa_1$, $\alpha$ and $\beta_1$.}
	\label{RQ_kappa1_alpha_beta1}
\end{figure}

Since RQ measures the tumor rings, we see in Figs $\ref{Ring_dif_kappa1}$-$\ref{Ring_dif_beta}$ how the model captures three kinds of tumor ring changing mainly the parameters $\kappa_1$, $\alpha$ or $\beta_1$. 
\\

Comparing Figs $\ref{Ring_dif_kappa1}$-$\ref{Ring_dif_beta}$, we appreciate as the change in $\alpha$ has a greater influence in tumor rings than $\kappa_1$ and $\beta$. Hence, the best configurations to obtain a slim (resp. thick) ring would be choose a big (resp. small) $\alpha$.% whereas if the ring is thicker, the configuration would be a small $\alpha$. 

\subsection{Density tumor growth}

Now, we measure the amount of density of total tumor in order to obtain the different tumor growths related to different values of $\kappa_1$, $\alpha$ and $\beta$.
\\

In Figure $\ref{densidad_razon_necrosis_tumor}$, we compare $\displaystyle\int_{\Omega}\left(T+N\right)\;dx$ for the different values of $\kappa_1$, $\alpha$ and $\beta_1$ chosen in Table $\ref{parametros_variables_razon_tumor_necrosis}$.
\begin{figure}[H]
	\centering
	\begin{subfigure}[b]{0.45\linewidth}
			\includegraphics[width=9cm, height=5cm]{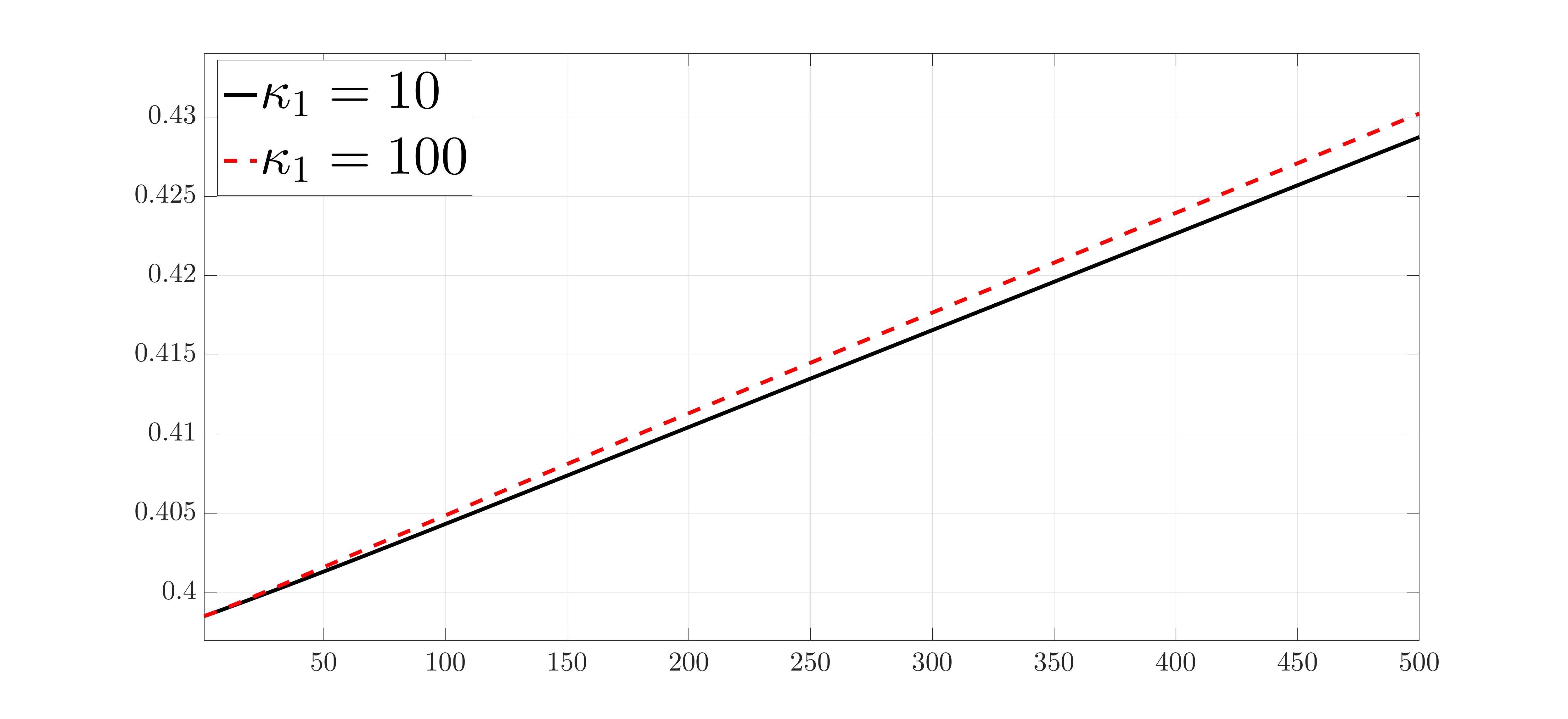}
		\centering
		\caption{$\displaystyle\int_{\Omega}\left(T+N\right)\;dx$ versus time for $\kappa_1$.}
		\label{densidad_anillo_kappa1}
	\end{subfigure}
	\hspace{1cm}
	\begin{subfigure}[b]{0.45\linewidth}
\includegraphics[width=9cm,height=5cm]{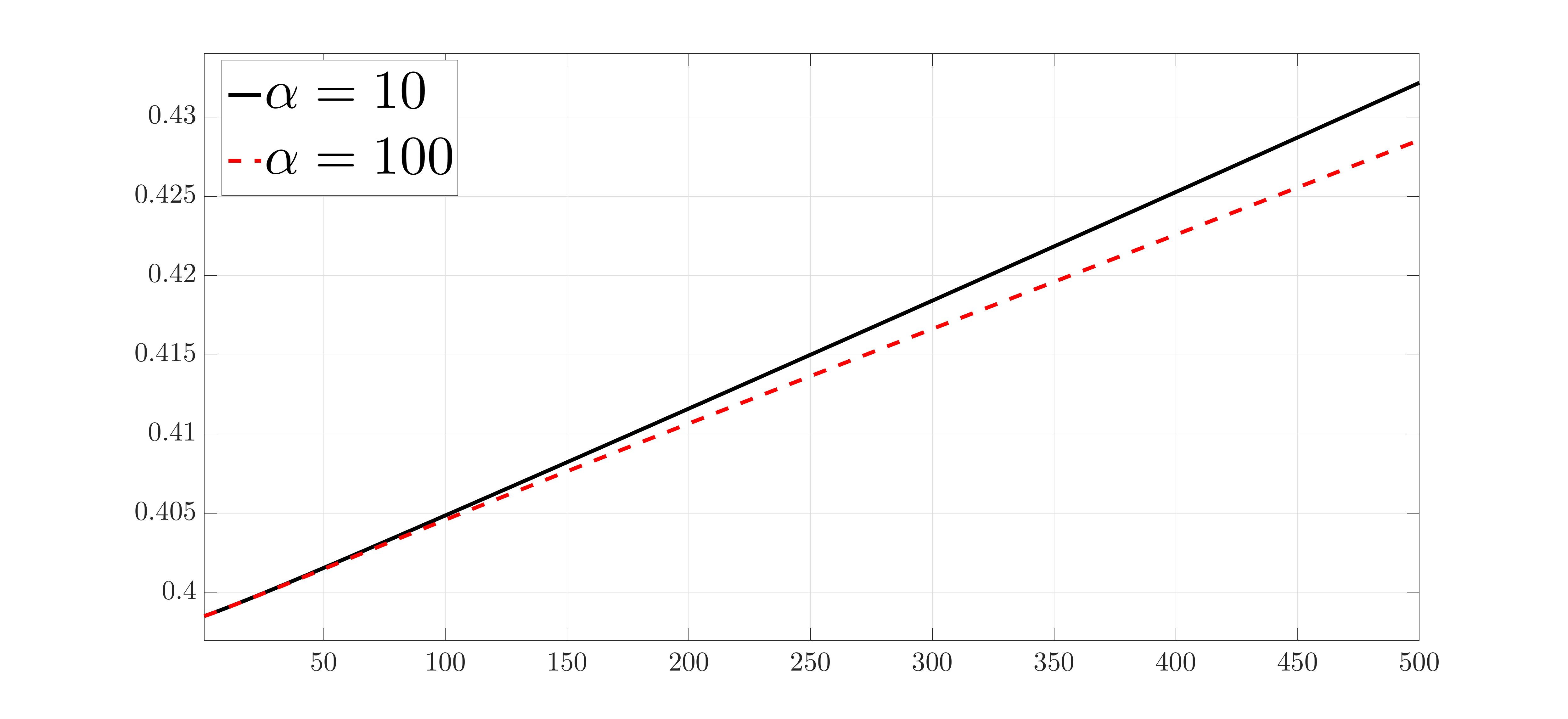}
		\centering
		\caption{$\displaystyle\int_{\Omega}\left(T+N\right)\;dx$ versus time for $\alpha$.}
		\label{densidad_anillo_alpha}
	\end{subfigure}
\end{figure}
%\vspace{15cm}
\begin{figure}[H]\ContinuedFloat
\begin{subfigure}[b]{0.45\linewidth}
	\includegraphics[width=9cm, height=5cm]{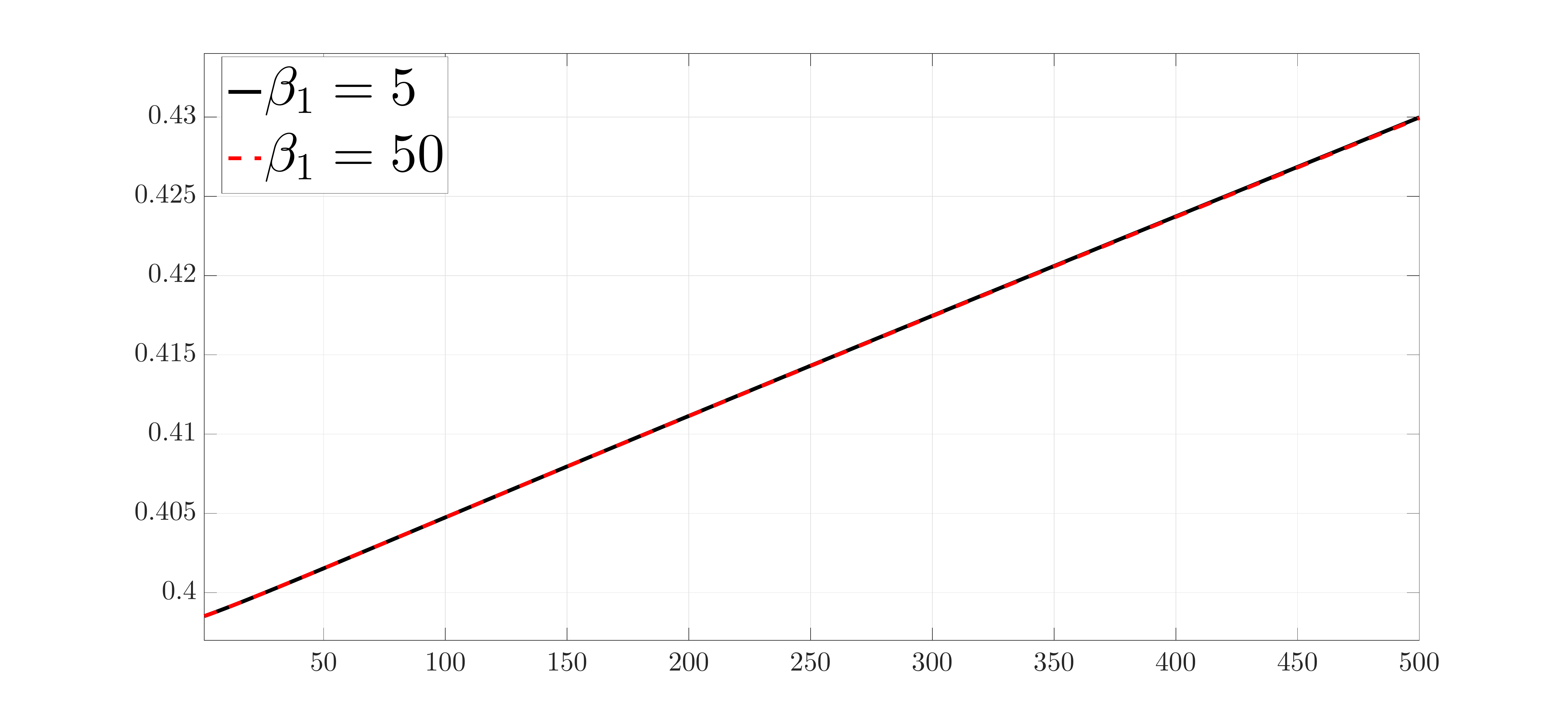}
	\centering
	\caption{$\displaystyle\int_{\Omega}\left(T+N\right)\;dx$ versus time for $\beta_1$.}
	\label{densidad_anillo_beta1}
\end{subfigure}
	\centering
	\caption{$\displaystyle\int_{\Omega}\left(T+N\right)\;dx$ versus time for $\kappa_1$, $\alpha$ and $\beta_1$.}
	\label{densidad_razon_necrosis_tumor}
\end{figure}
We see in Figure $\ref{densidad_razon_necrosis_tumor}$ that the parameters $\kappa_1$ and $\alpha$ produce more variation than $\beta_1$ in the total tumor density.
\\

With respect to $\alpha$, see Figure $\ref{densidad_anillo_alpha}$, we get the maximum density for $\alpha=10$ due to low hypoxia allows a higher tumor growth than for high hypoxia. Thus, the minimum density is obtained for $\alpha=100$. In the case of $\kappa_1$, see Figure $\ref{densidad_anillo_kappa1}$, the difference between the two total densities is lower than for $\alpha$. The maximum and minimum densities are achieved for $\kappa_1=100$ and $\kappa_1=10$, respectively.
\\

Furthermore, the highest density is obtained for $\alpha=10$ while the density for $\alpha=100$ and $\kappa_1=10$ are similar.

\subsection{Conclusions}

Based on the study \cite{Julian_2016}, we know that tumors with a thick tumor ring have the worst prognosis. This kind of tumor growth in $\left(\ref{probOriginal}\right)$-$\left(\ref{condinicio}\right)$ is produced by a low value of the parameter $\alpha$. This means that a change of the rate of tumor destruction for hypoxia produces much difference in the tumor rings, as we see in Figure $\ref{Ring_dif_alpha}$. In addition, we have obtained more total density in the case of low $\alpha$ as we see in Figure $\ref{densidad_anillo_alpha}$.
\\

With respect to $\kappa_1$, there also exist differences in the tumor rings and total density for different values of $\kappa_1$ observed in Figs $\ref{Ring_dif_kappa1}$ and $\ref{densidad_anillo_kappa1}$. However, these differences are not as relevant as the changes with variations for $\alpha$. 
\\

The parameter related to the tumor destruction by necrosis, $\beta_1$, is not relevant in this prognosis since $\beta_1$ does not distinguish tumors with thick and slim tumor rings and the variation of $\beta_1$ does not produce changes in the total density as we can see in Figure $\ref{densidad_anillo_beta1}$.
\\

In conclusion, $\alpha$ is the most important parameter in relation to the tumor ring and $\alpha$ and $\kappa_1$ have more relevancy for total density in the tumor growth than $\beta_1$.
\section{Regularity surface}\label{regularidad}

Our model $\left(\ref{prob_K_rho}\right)$ can get different regularities for the tumor surfaces as we see in Figure $\ref{Sim_crec_reg_irreg}$, where we show in Figure $\ref{crec_redondo}$ the tumor growth with the initial vasculature uniformly distributed in the space, similar to the previous section, and in Figure $\ref{crec_flor}$ the tumor growth with the initial vasculature distributed in three zones with different concentrations of vasculature:

\begin{figure}[H]
	\centering
	\begin{subfigure}[b]{0.45\linewidth}
				\includegraphics[width=\linewidth]{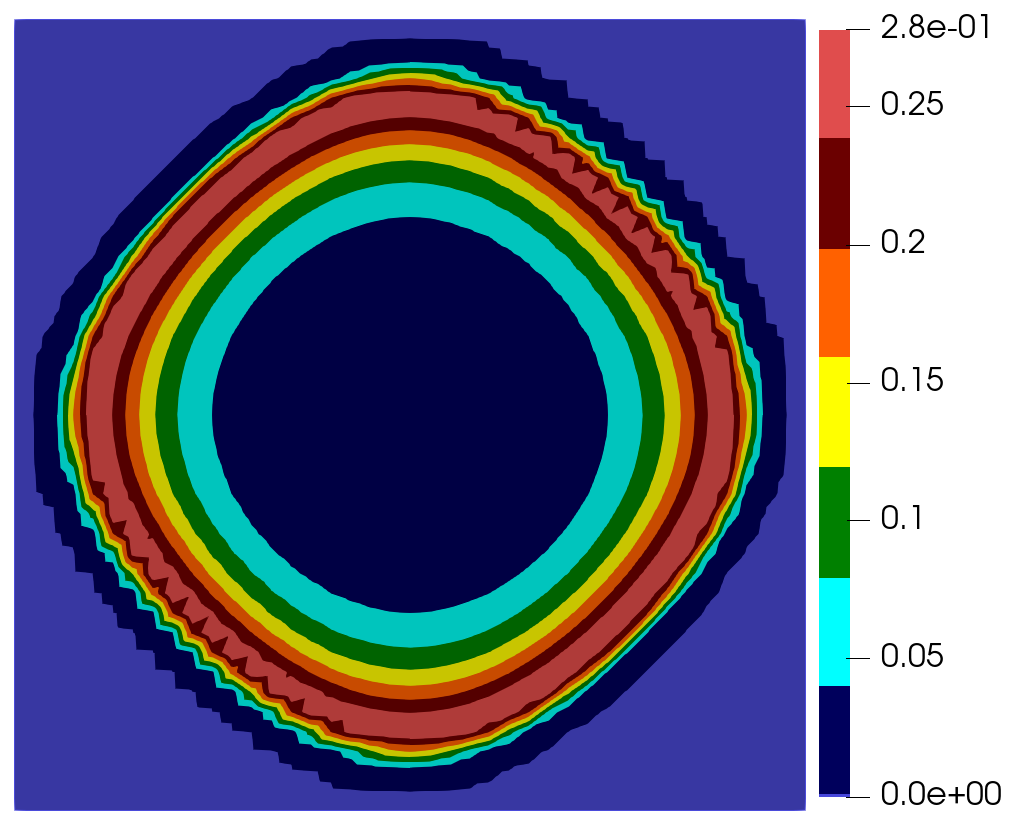}
		\centering
		\caption{Tumor growth with vasculature \\
			\centering uniformly distributed.}
		\label{crec_redondo}
	\end{subfigure}
	\hspace{1cm}
	\begin{subfigure}[b]{0.45\linewidth}
			\includegraphics[width=\linewidth]{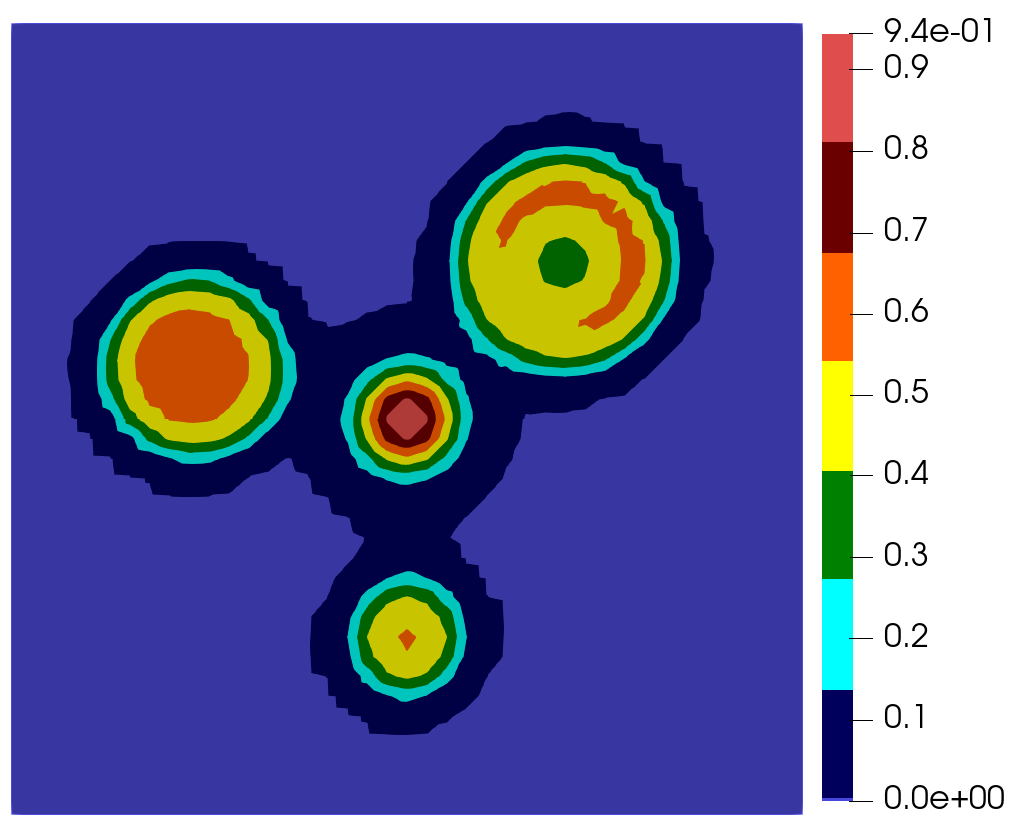}
		\centering
		\caption{Tumor growth with vasculature \\
			\centering non-uniformly distributed.}
		\label{crec_flor}
	\end{subfigure}
	\caption{Tumor growth for different distributions of vasculature}
	\label{Sim_crec_reg_irreg}
\end{figure}

Thus, we see as our model captures tumor growths with regular or irregular surface due to the speed of the tumor diffusion function of tumor depends on the vasculature. We remember that the adimensionalized diffusion term is defined by $\nabla\cdot\left(\left(\kappa_1\;P\left(\Phi,T\right)+1\right)\nabla\;T\right)$ with $\kappa_1>0$ and $P\left(\Phi,T\right)$ defined in $\left(\ref{funcionP}\right)$. In particular, the parameter $\kappa_1$ regulates the influence of the vasculature spatial distribution in the  regularity surface of the tumor.  
\\

In order to detect which parameter could be more effective for the regularity surface of the tumor, we show some simulations moving the value of one parameter and checking how the tumor growth changes.
\\

Since we focus our criteria on tumor and vasculature, we are going to move the parameters which appear in their equations, that is, $\kappa_1$, $\alpha$, $\beta_1$, $\beta_2$, $\gamma$ and $\delta$. 
\\

For these parameters we take the following values:

\begin{table}[H]
	\centering
	\begin{tabular}{c|c|c|c|c|c|c}
		%\hline
		\textbf{Variable (Fixed value)} &	$\kappa_1\;\;\left(55\right)$  & 	$\alpha\;\;\left(45\right)$  & 	$\beta_1\;\;\left(27.5\right)$  &	$\beta_2\;\;\left(2.55\right)$  &$\gamma\;\;\left(0.255\right)$  &$\delta\;\;\left(2.55\right)$   \\
		\hline
		\textbf{Ranges}  & $\left[10,\;100\right]$ & $\left[10,\;100\right]$ & $\left[5,\;50\right]$ & $\left[0.1,\;5\right]$ & 	$\left[0.01,\;0.5\right]$ & $\left[0.1,\;5\right]$  
	\end{tabular}
	\caption{\label{parametros_variables_crec_reg_irreg} Variable value parameters.}
\end{table}

In the following simulations, the initial vasculature is distributed in various zones with different concentrations along the domain as we can see in Figure $\ref{vasc_incial}$:
\begin{figure}[H]
	\includegraphics[width=8cm, height=5cm]{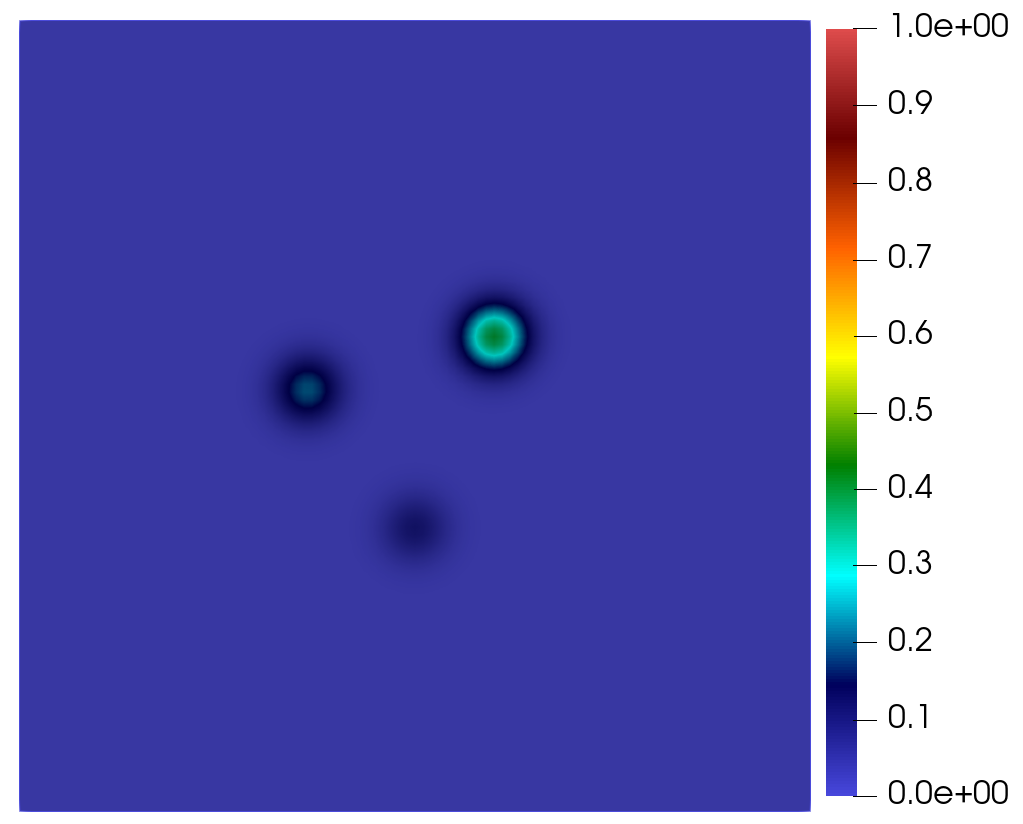}
	\centering
	\caption{Initial vasculature.}
	\label{vasc_incial}
\end{figure}

We remember that initial tumor is defined as in Figure $\ref{tumor_incial}$ and necrosis is initially zero.
\subsection{Regularity Surface quotient}
We start showing the simulations according to the following quotient between the area occupied by the total tumor (tumor and necrosis) and the area of a sphere whose radio is equal to the maximum radio of the tumor, that is the smallest sphere containing the tumor. Thus, we show this difference for the different values of $\kappa_1$, $\alpha$, $\beta_1$, $\gamma$, $\delta$ and $\beta_2$ chosen in Table $\ref{parametros_variables_crec_reg_irreg}$. For this, we have considered the following "surface quotient" (SQ) coefficient:

\begin{equation}\label{SQ}
0\leq\text{SQ}=\dfrac{\displaystyle\int_{\Omega}\left(T+N\right)_{\min}\;dx}{\pi\cdot \;\left(\textbf{R}_{\max}\right)^2}\leq 1
\end{equation}
where $\left(T+N\right)_{\min}$ and $\textbf{R}_{\max}$ are defined as follows: 
\begin{equation}\label{Tmin}
\left(T+N\right)_{\min}=\left\{
\begin{array}{ll}
1&\text{if}\;\;T+N\geq 0.001,\\
\\
0&\text{otherwise}. 
\end{array}
\right.
\end{equation}

\begin{equation}\label{Rmax}
\textbf{R}_{\max}=\max\left\{\text{radio of the subdomain where }\left(T+N\right)_{\min}=1\right\}.
\end{equation}
%\begin{figure}[H]
%	\includegraphics[width=18cm, height=20cm]{CREC_TUMOR/EPS/Crec_Reg_Irreg_kappa1_rho_alpha_beta_delta_gamma.eps}
%	\centering
%	\caption{Difference between volumen of total tumor and a sphere of maximum radio for different values of $\kappa_1$, $\rho$, $\alpha$, $\beta$, $\gamma$ and $\delta$}
%	\label{crec_irreg_reg}
%\end{figure}

Thus, we could conclude that if SQ is near to zero, the surface will be irregular whereas if SQ is close to one, the surface will be regular.

\begin{figure}[H]
	\centering
	\begin{subfigure}[b]{0.45\linewidth}
		\includegraphics[width=9cm, height=5cm]{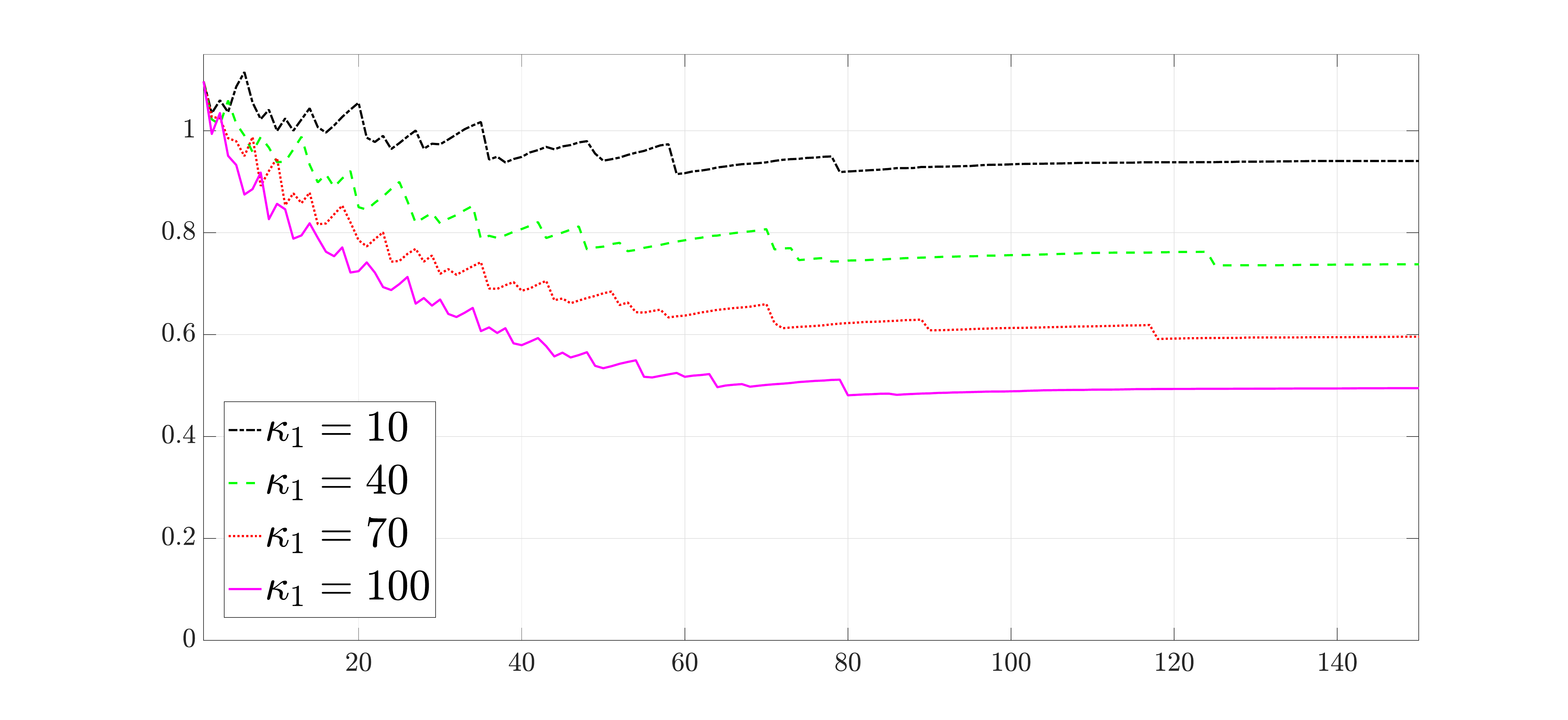}
		\centering
		\caption{SQ versus time for $\kappa_1$.}
		\label{crec_irreg_kappa1}
	\end{subfigure}
	\hspace{1cm}
	\begin{subfigure}[b]{0.45\linewidth}
		\includegraphics[width=9cm, height=5cm]{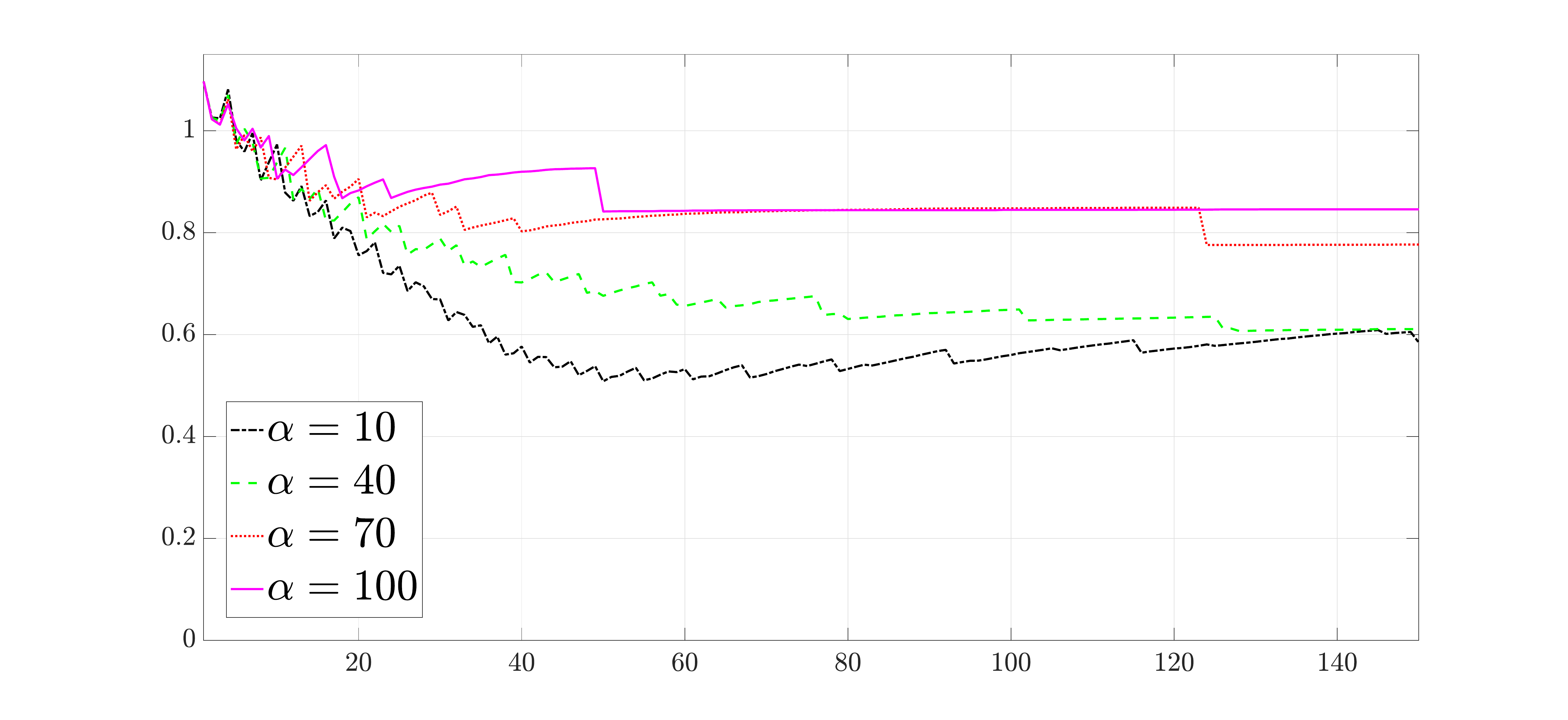}
		\centering
		\caption{SQ versus time for $\alpha$.}
		\label{crec_irreg_alpha}
	\end{subfigure}
\end{figure}
%\vspace{15cm}
\begin{figure}[H]\ContinuedFloat
	\begin{subfigure}[b]{0.45\linewidth}
			\includegraphics[width=9cm, height=5cm]{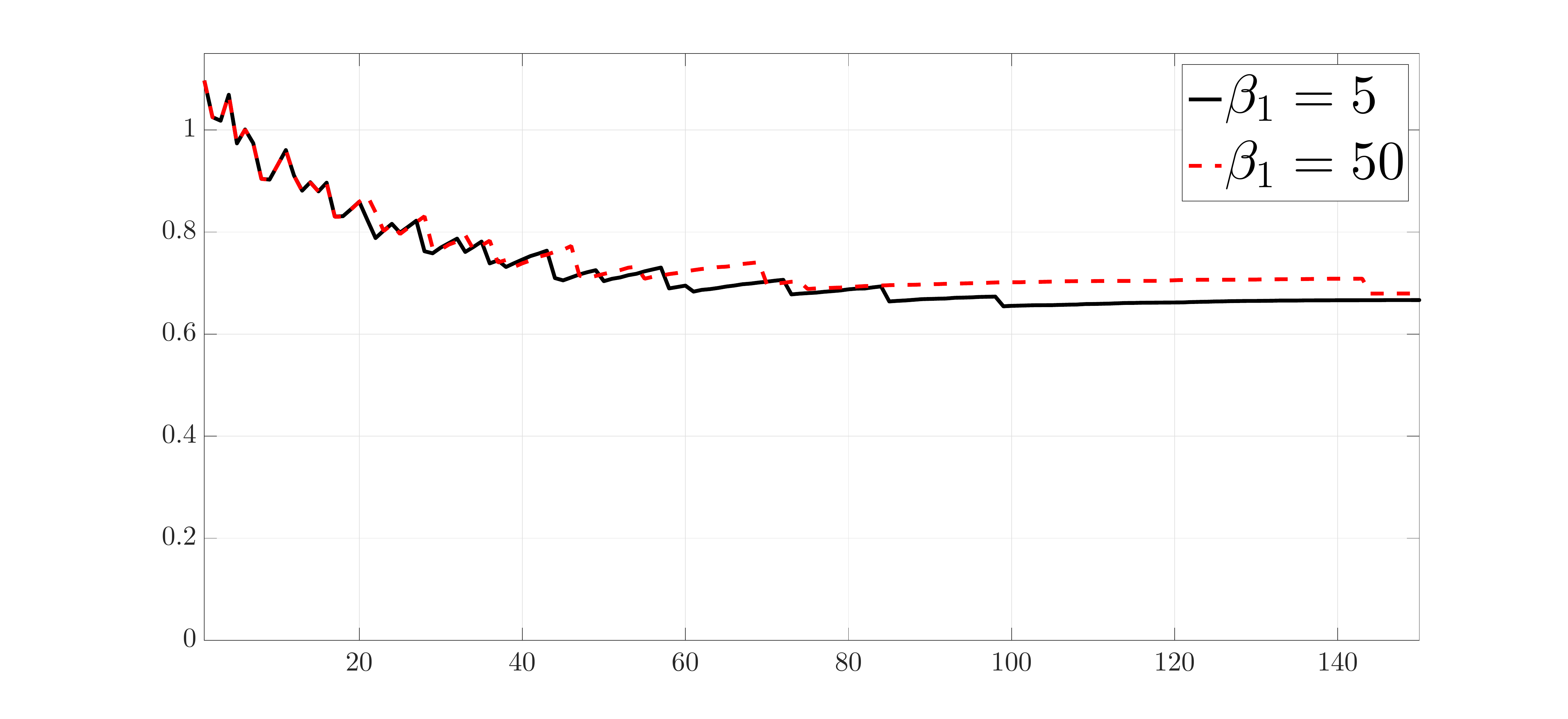}
		\centering
		\caption{SQ versus time for $\beta_1$.}
		\label{crec_irreg_beta1}
	\end{subfigure}
	\hspace{1cm}
	\begin{subfigure}[b]{0.45\linewidth}
	\includegraphics[width=9cm, height=5cm]{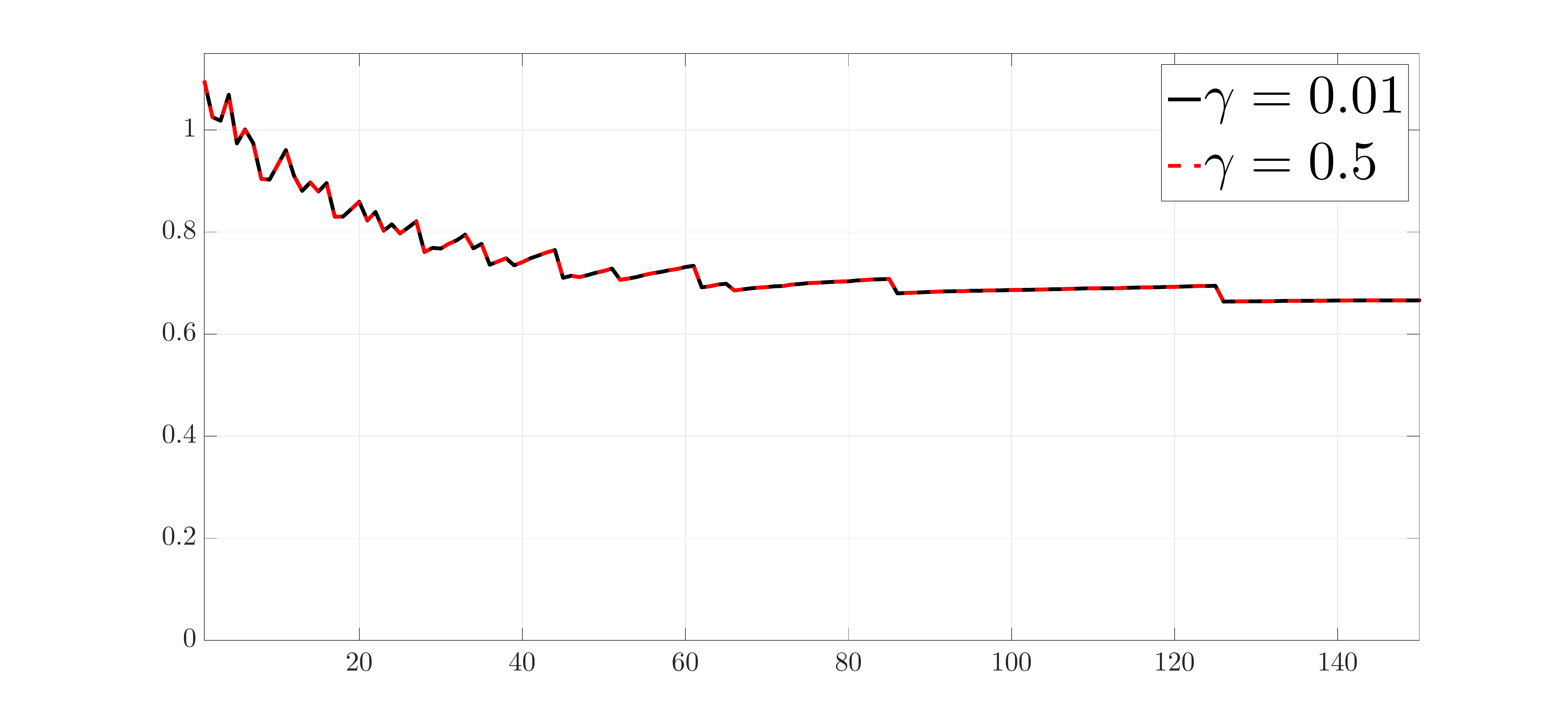}
	\centering
	\caption{SQ versus time for $\gamma$.}
	\label{crec_irreg_gamma}
\end{subfigure}
\end{figure}
\begin{figure}[H]\ContinuedFloat
\begin{subfigure}[b]{0.45\linewidth}
	\includegraphics[width=9cm, height=5cm]{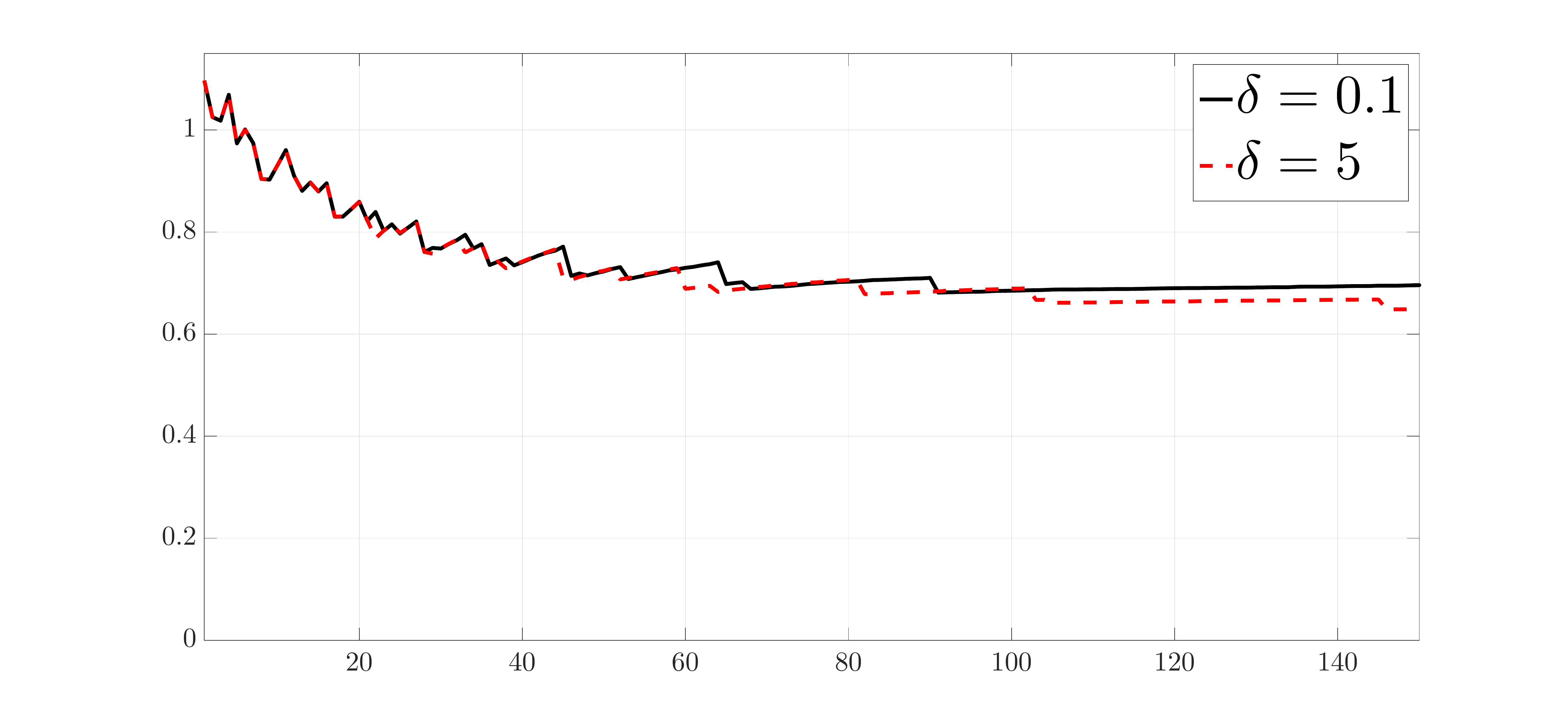}
	\centering
	\caption{SQ versus time for $\delta$.}
	\label{crec_irreg_delta}
\end{subfigure}
	\hspace{1cm}
\begin{subfigure}[b]{0.45\linewidth}
	\includegraphics[width=9cm, height=5cm]{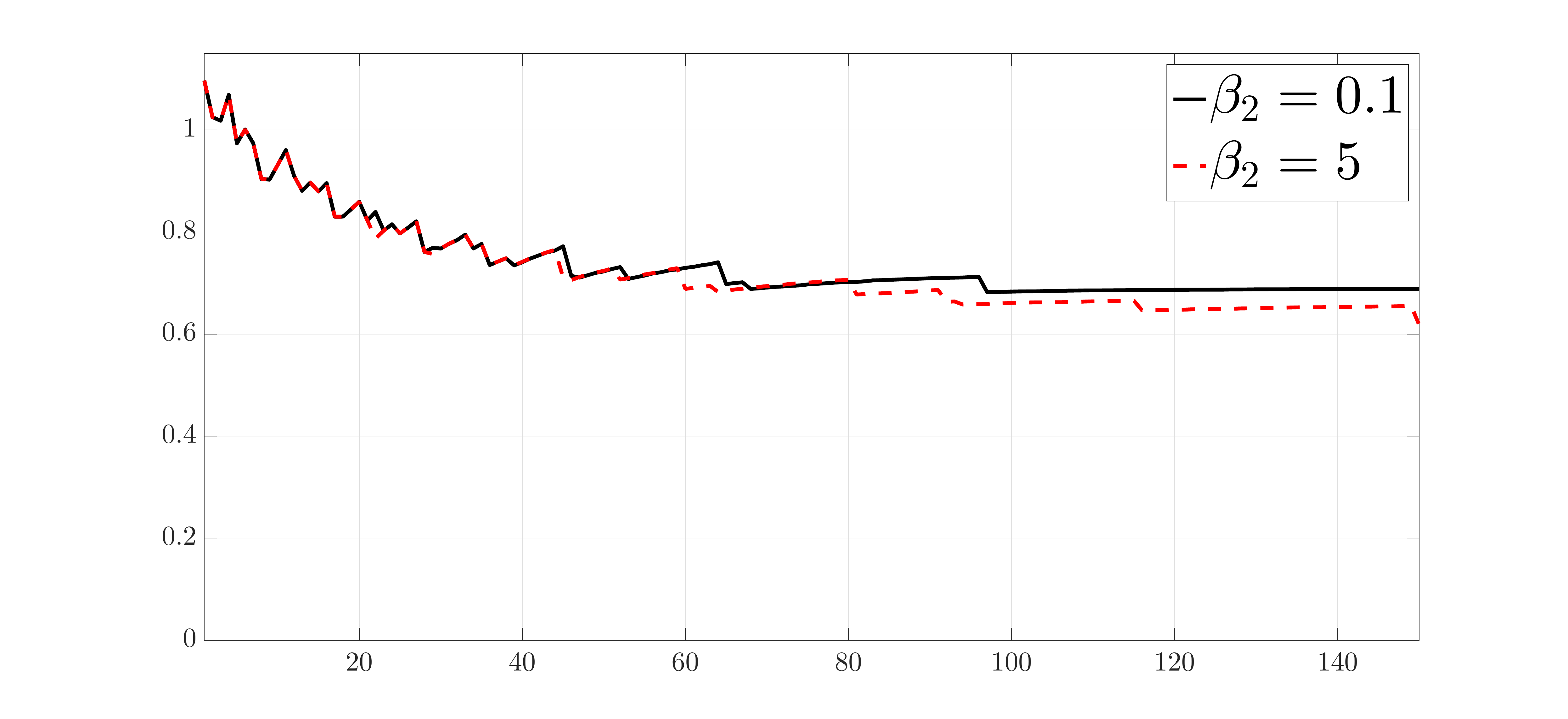}
	\centering
	\caption{SQ versus time for $\beta_2$.}
	\label{crec_irreg_beta2}
\end{subfigure}
	\centering
	\caption[]{SQ versus time for $\kappa_1$, $\alpha$, $\beta_1$, $\gamma$, $\delta$ and $\beta_2$.}
	\label{SQ_kappa1_alpha_beta1_gamma_delta_beta2_2}
\end{figure}
\begin{obs}
	Due to the size of mesh considered, at the beginning of the pictures given in Figure $\ref{SQ_kappa1_alpha_beta1_gamma_delta_beta2_2}$, the value of SQ is larger than $1$ and it is observed oscillations in the graphs of SQ. Indeed, if we consider a mesh size smaller, these initial values of SQ and the oscillations can be corrected. In order to demonstrate this, we show an example of SQ versus time considering a mesh size smaller:
	\begin{figure}[H]
		\centering
		\includegraphics[width=9cm, height=5cm]{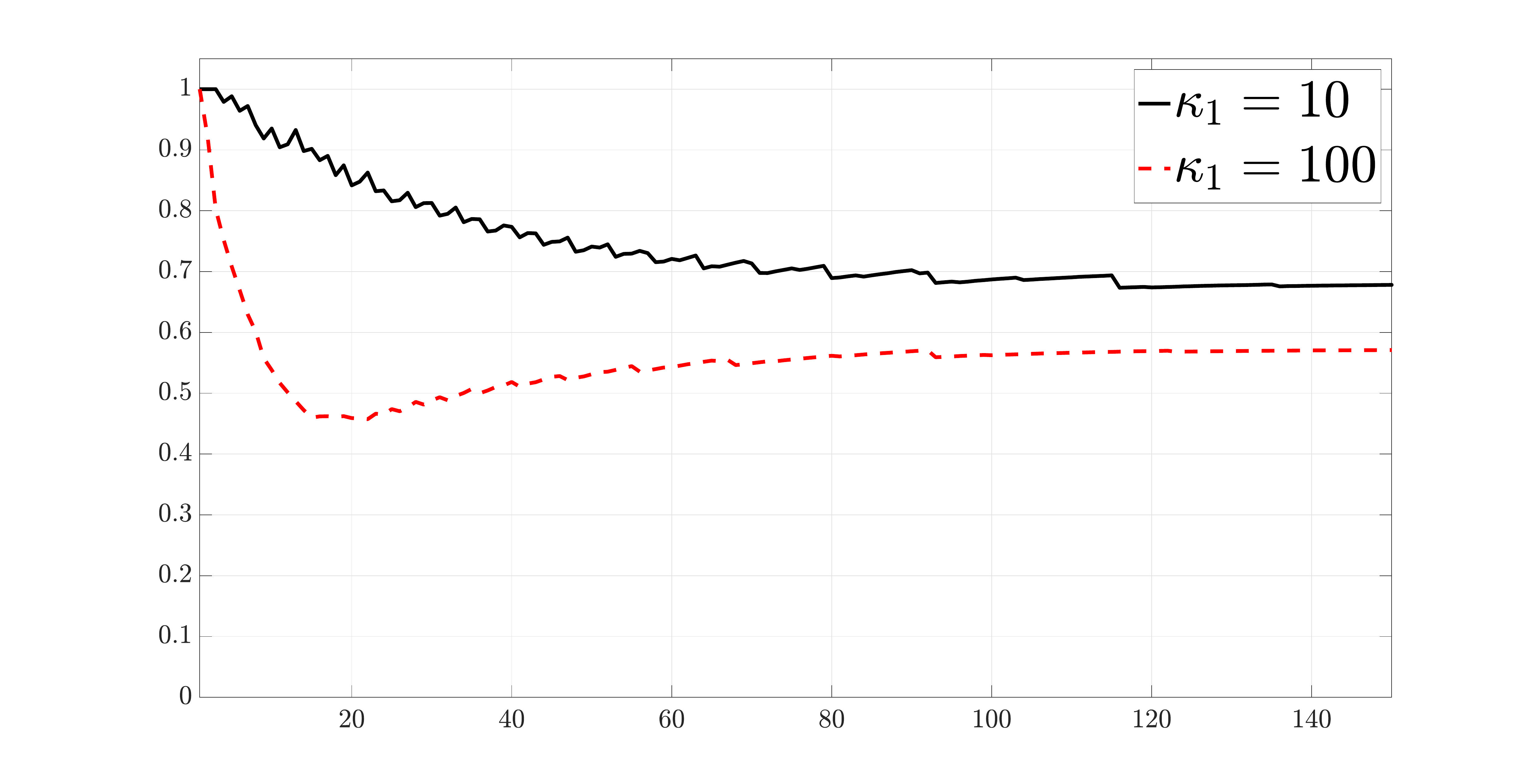}
		\centering
		\caption{SQ versus time for $\kappa_1$ for a mesh size smaller.}
		\label{crec_irreg_kappa_malla_fina}
	\end{figure}
	
	However, it is not completely necessary the use of a mesh size smaller since we obtain the same behaviour (in average) for the mesh considered initially and we reduce the computational time. 
\end{obs}
We see in Figs $\ref{crec_irreg_kappa1}$-$\ref{crec_irreg_beta2}$ how our model can differentiate two kinds of tumor growth changing mainly the parameters $\kappa_1$ and/or $\alpha$, see Figs $\ref{crec_irreg_kappa1}$, and $\ref{crec_irreg_alpha}$, respectively while the variation in the parameters $\beta_1$, $\gamma$, $\delta$ and $\beta_2$ do not change the irregularity of tumor growth as we see in Figs $\ref{crec_irreg_beta1}$-$\ref{crec_irreg_beta2}$.
\\

Once we have identified that the more important parameters for the regularity surface are $\kappa_1$ and $\alpha$, we measure the area of total tumor for the different values of $\kappa_1$ and $\alpha$ taken in Figs $\ref{crec_irreg_kappa1}$ and $\ref{crec_irreg_alpha}$.

\subsection{Tumor area}

Here, we will compare the area occupied by total tumor for the different values of $\kappa_1$ and $\alpha$. In order to measure this area, we consider:

\begin{equation}\label{Area}
\displaystyle\int_{\Omega}\left(T+N\right)_{\min}\;dx
\end{equation}
where $\left(T+N\right)_{\min}$ is defined in $\left(\ref{Tmin}\right)$.
\\

In the following graphs we show the area of total tumor (tumor and necrosis) for the different values of $\kappa_1$ and $\alpha$ chosen in Table $\ref{parametros_variables_crec_reg_irreg}$. 

\begin{figure}[H]
\centering
\begin{subfigure}[b]{0.45\linewidth}
	\includegraphics[width=9cm, height=5cm]{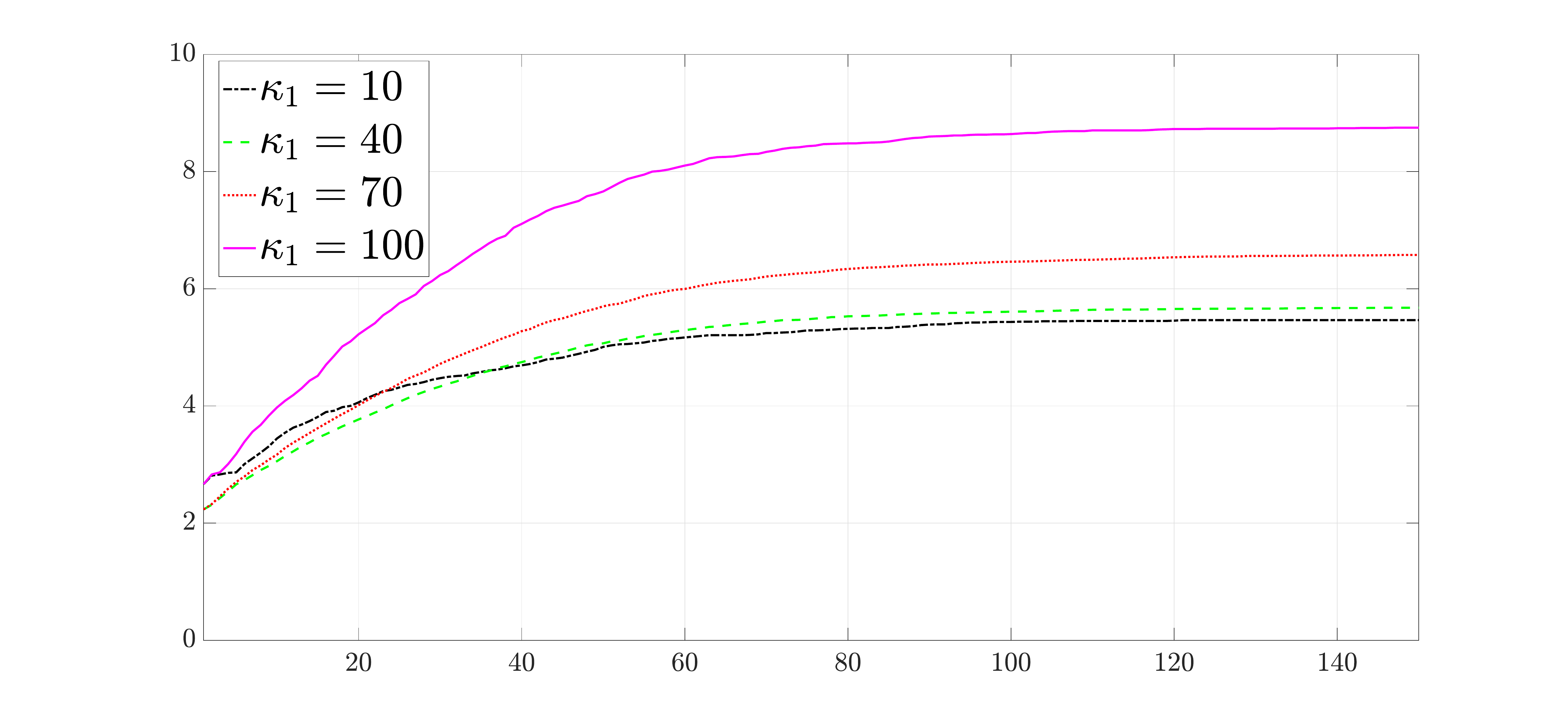}
	\centering
	\caption{Area of total tumor versus time for $\kappa_1$.}
	\label{area_kappa1}
\end{subfigure}
\hspace{1cm}
\begin{subfigure}[b]{0.45\linewidth}
	\includegraphics[width=9cm, height=5cm]{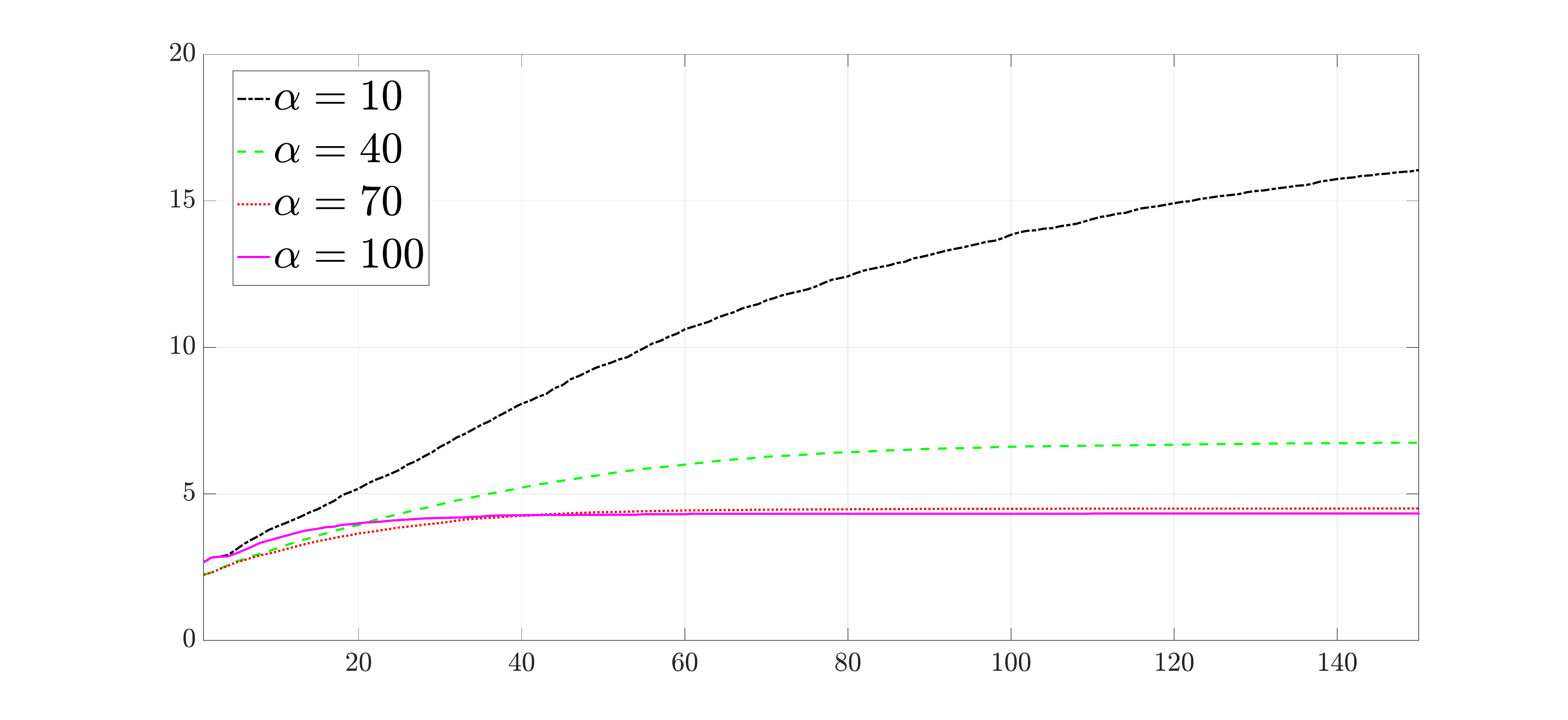}
	\centering
	\caption{Area of total tumor time for $\alpha$.}
	\label{area_alpha}
\end{subfigure}
	\centering
	\caption{Area of total tumor versus time for $\kappa_1$ and $\alpha$.}
	\label{area_crec_necrosis_tumor}
\end{figure}

We see in Figure $\ref{area_crec_necrosis_tumor}$ how the largest area correspond to $\alpha=10$. For $\kappa_1=100$, we also obtain a high area due to the large value of anisotropic speed diffusion $\kappa_1$. Finally, for $\kappa_1=10$ and $\alpha=100$, total area has a similar variation due to the effect of low diffusion in the case of $\kappa_1=10$ and high tumor destruction for hypoxia in the case of $\alpha=100$.
\\

Thus, we deduce that the variation speed of total area is not constant for the different values of $\kappa_1$ and $\alpha$ considered in Figure $\ref{area_crec_necrosis_tumor}$. However, if we consider the "surface quotient" (SQ), we obtain a similar variation between the different values of $\kappa_1$ and $\alpha$, see Figure $\ref{SQ_kappa1_alpha_beta1_gamma_delta_beta2_2}$. Hence, the factor which modifies this variation is the term $\textbf{R}_{\max}$, defined by $\ref{Rmax}$ . Moreover, $\textbf{R}_{\max}$ will change more with the variation of $\alpha$ that for different values of $\kappa_1$ since the variation of SQ for $\kappa_1$, Figure $\ref{crec_irreg_kappa1}$, is bigger  than for $\alpha$, Figure $\ref{crec_irreg_alpha}$.

\subsection{Tumor growth}
In this part, we will show the tumor growth for $\kappa_1=100$ and $\alpha=10$ in five time steps in order to see the spatial growth of tumor. The rest of parameters take the values showed in Table $\ref{parametros_variables_crec_reg_irreg}$.
\\

\begin{figure}[H]
	\centering
	\begin{subfigure}[b]{0.15\linewidth}
	\includegraphics[width=1.2\linewidth]{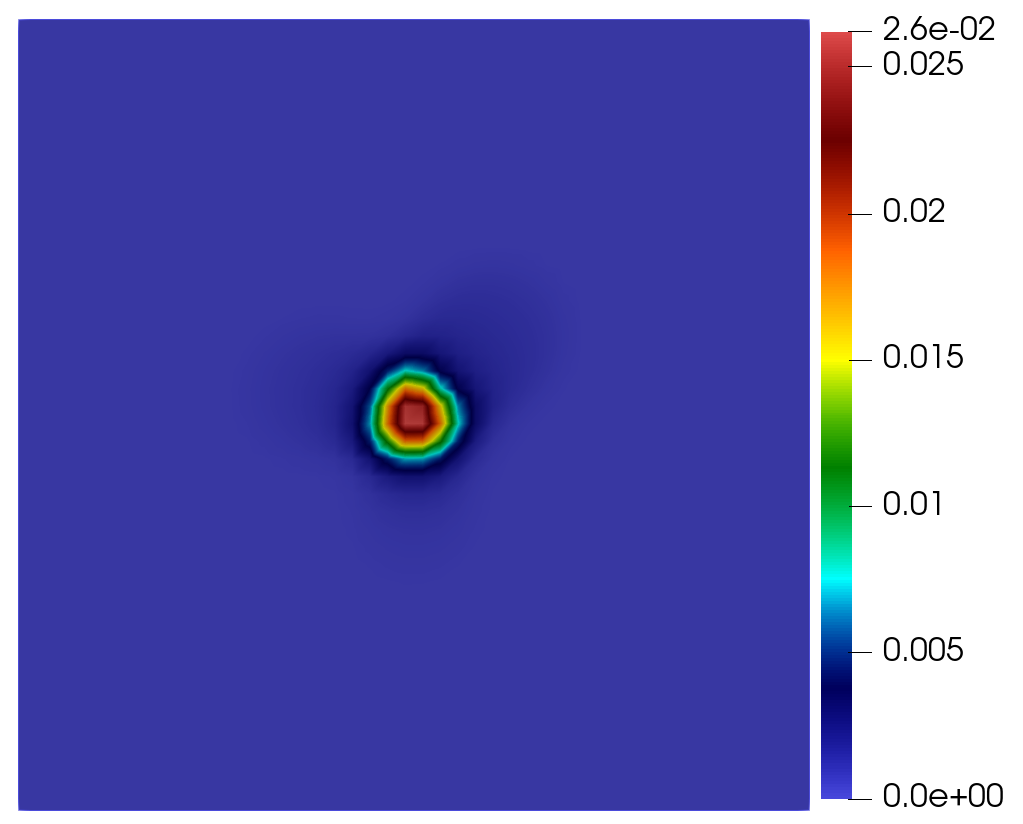}
		\centering
		\caption{$t=50$}
	\end{subfigure}
	\hspace{0.4cm}
	\begin{subfigure}[b]{0.15\linewidth}
		\includegraphics[width=1.2\linewidth]{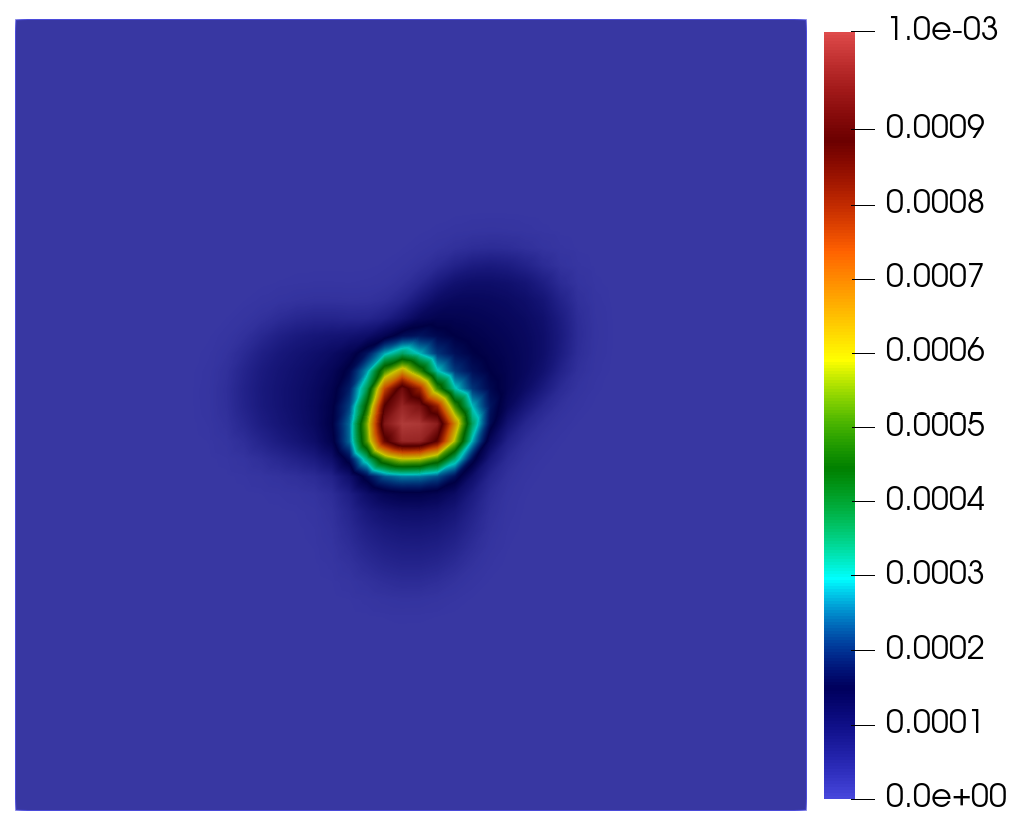}
		\centering
		\caption{$t=100$}
		\label{kappa1_100_b}
	\end{subfigure}
	\hspace{0.4cm}
	\begin{subfigure}[b]{0.15\linewidth}
	\includegraphics[width=1.2\linewidth]{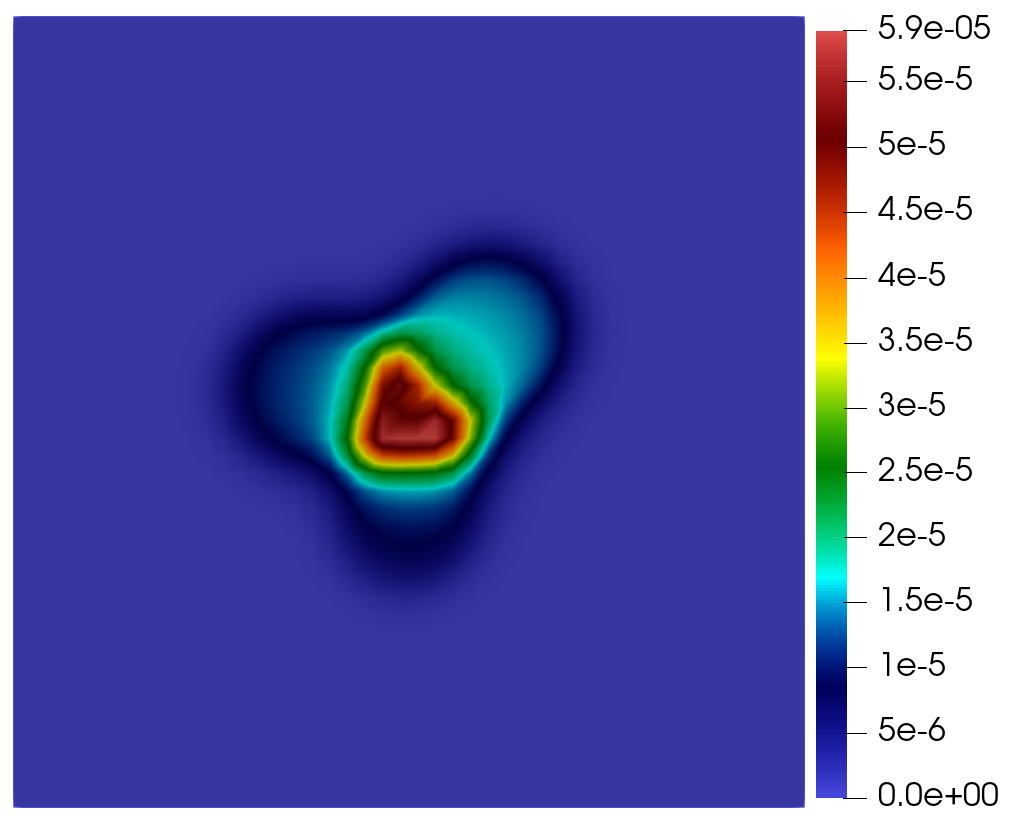}
		\centering
		\caption{$t=150$}
	\end{subfigure}
	\hspace{0.4cm}
	\begin{subfigure}[b]{0.15\linewidth}
		\includegraphics[width=1.2\linewidth]{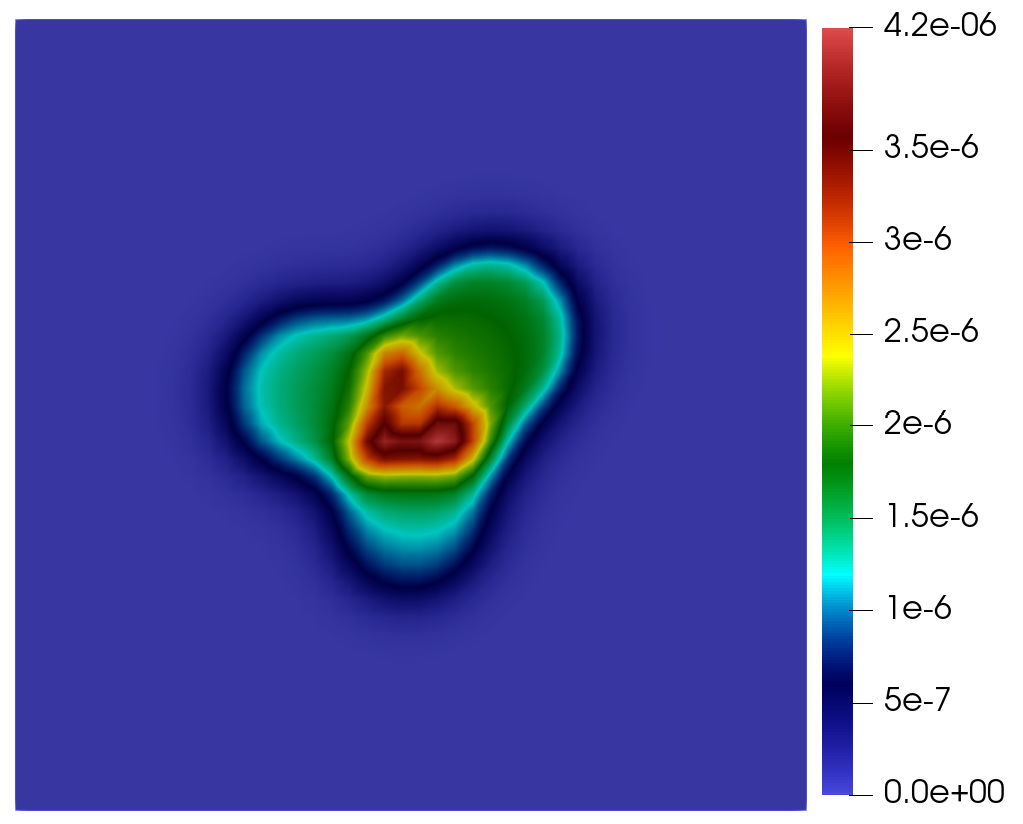}
		\centering
		\caption{$t=200$}
	\end{subfigure}
	\hspace{0.4cm}
	\begin{subfigure}[b]{0.15\linewidth}
		\includegraphics[width=1.2\linewidth]{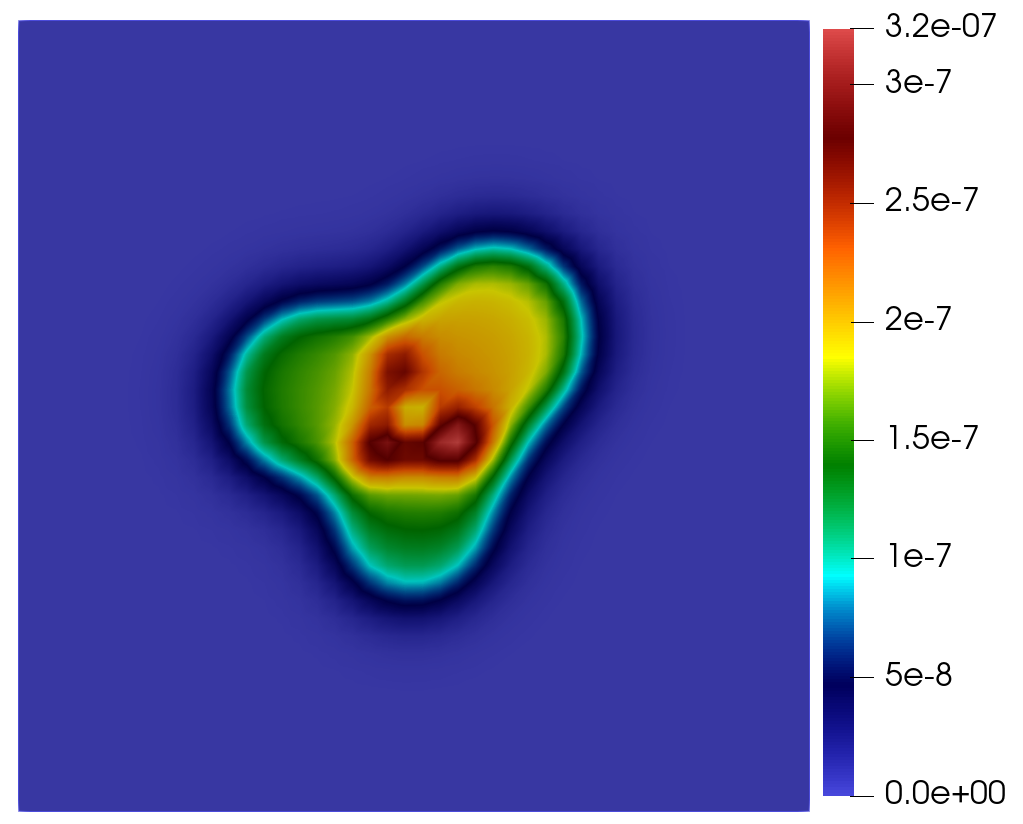}
		\centering
		\caption{$t=250$}
		\label{kappa1_100_e}
	\end{subfigure}
	\caption{Irregular tumor growth for $\kappa_1=100$.}
	\label{crec_irreg_kappa1_100}
\end{figure}
\begin{figure}[H]
	\centering
	\begin{subfigure}[b]{0.15\linewidth}
		\includegraphics[width=1.2\linewidth]{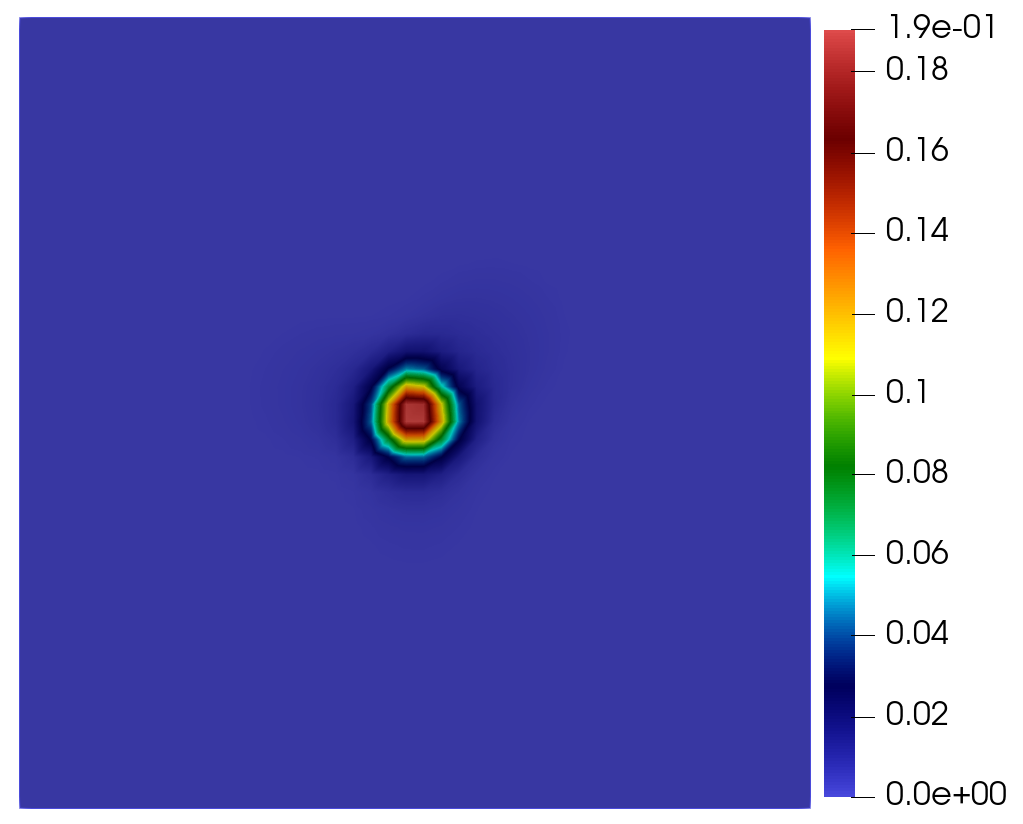}
		\centering
		\caption{$t=50$}
	\end{subfigure}
	\hspace{0.4cm}
	\begin{subfigure}[b]{0.15\linewidth}
		\includegraphics[width=1.2\linewidth]{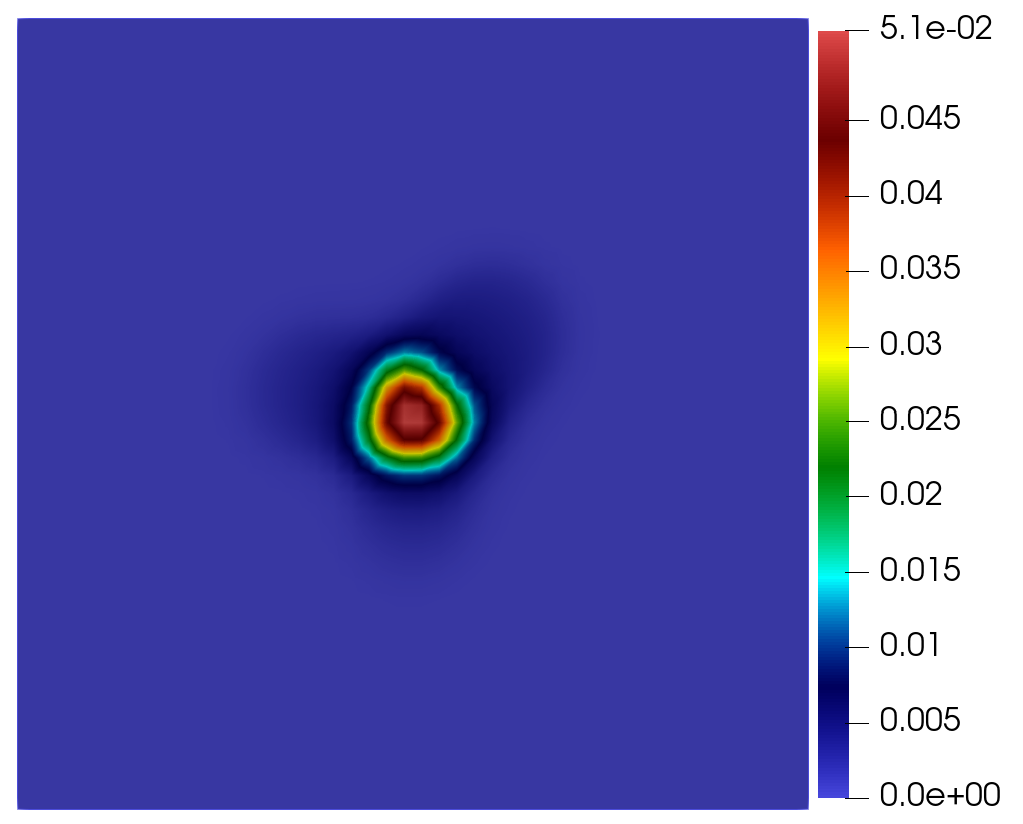}
		\centering
		\caption{$t=100$}
	\end{subfigure}
	\hspace{0.4cm}
	\begin{subfigure}[b]{0.15\linewidth}
	\includegraphics[width=1.2\linewidth]{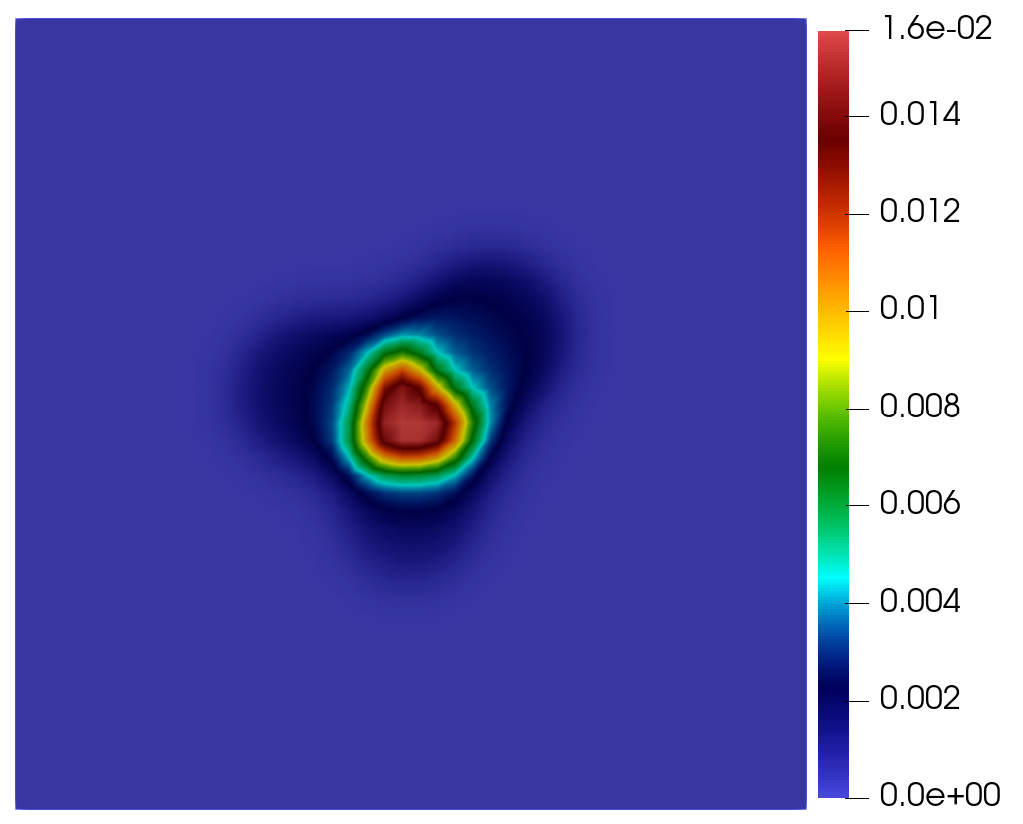}
		\centering
		\caption{$t=150$}
	\end{subfigure}
	\hspace{0.4cm}
	\begin{subfigure}[b]{0.15\linewidth}
	\includegraphics[width=1.2\linewidth]{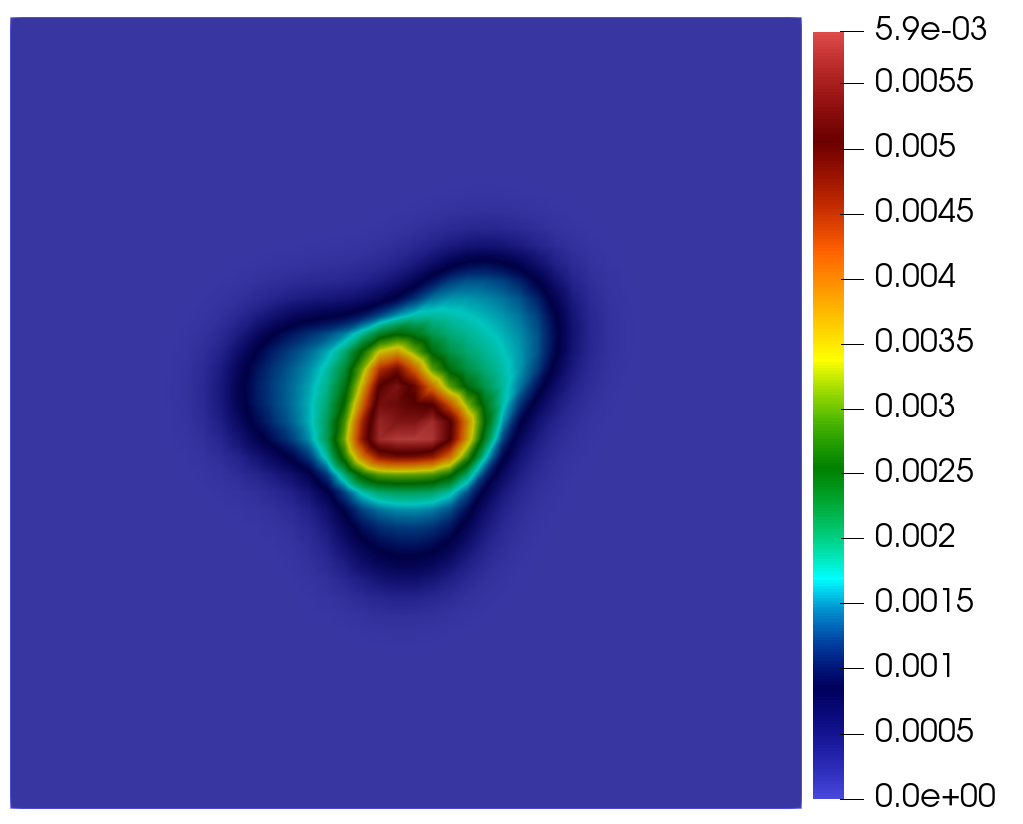}
		\centering
		\caption{$t=200$}
	\end{subfigure}
	\hspace{0.4cm}
	\begin{subfigure}[b]{0.15\linewidth}
	\includegraphics[width=1.2\linewidth]{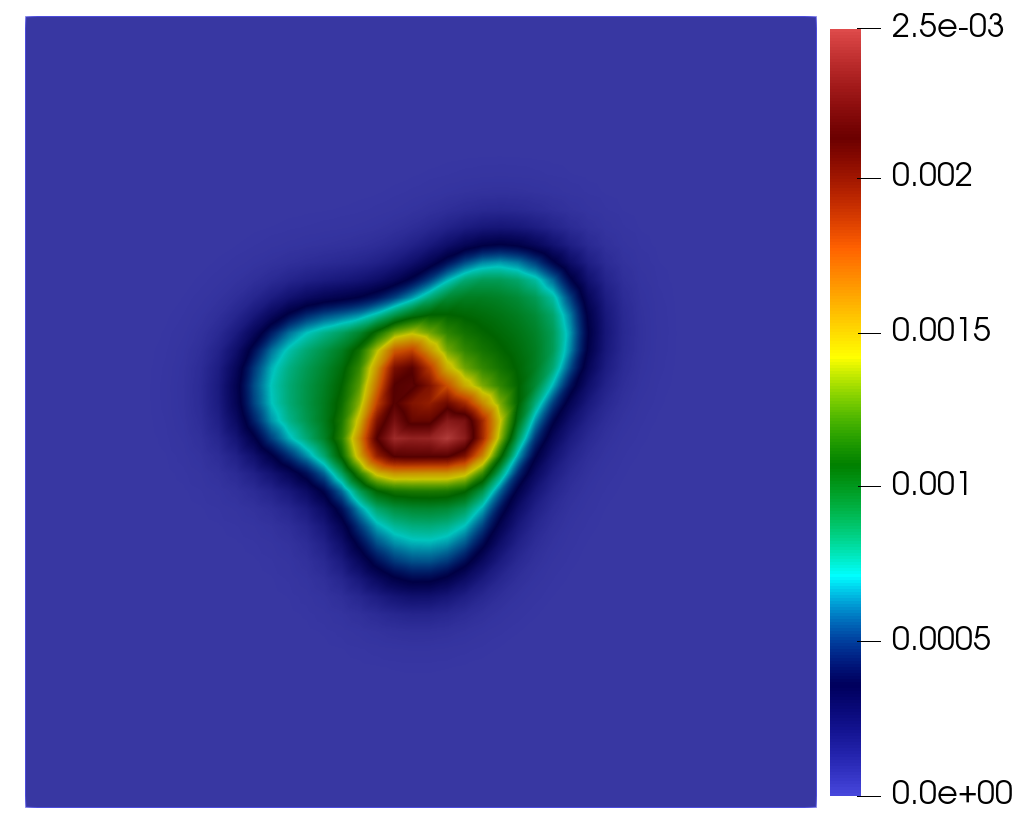}
		\centering
		\caption{$t=250$}
		\label{alpha_10_e}
	\end{subfigure}
	\caption{Tumor growth for $\alpha=10$.}
	\label{crec_irreg_alpha_10}
\end{figure}

We observe a faster tumor growth for $\kappa_1=100$ than for $\alpha=10$ in each time step. However, the amount of tumor for $\alpha=10$ is higher than for $\kappa_1=100$. These results are in concordance with the obtained in Figs $\ref{crec_irreg_kappa1}$ and $\ref{crec_irreg_alpha}$, where we see more irregularity for $\kappa_1=100$ than for $\alpha=10$, and Figs $\ref{area_kappa1}$ and $\ref{area_alpha}$, where the area for $\alpha=10$ is higher than for $\kappa_1=100$. Furthermore, from Figure $\ref{kappa1_100_b}$, we observe that the maximum value of tumor is lower than the critical value of $0.001$ given in $\left(\ref{Tmin}\right)$ whereas in Figure $\ref{crec_irreg_alpha_10}$ there are always zones where tumor achieves this critical value. Hence, it is normal that the area of total tumor for $\alpha=10$ be higher than for $\kappa_1=100$.

\subsection{Conclusions}

Based on the study \cite{Victor_2018}, tumors with a high irregularity in their surface have the worst prognosis. This correspond to tumors with a large value of $\kappa_1$ as we see in Figure $\ref{crec_irreg_kappa1}$. However, a low value of $\alpha$ also produces irregularity in the tumor surface, see Figure $\ref{crec_irreg_alpha}$ and a higher amount of total area, see Figure  $\ref{area_alpha}$, than for large $\kappa_1$, see Figure  $\ref{area_kappa1}$.
\\

Finally, we conclude that $\kappa_1$ and $\alpha$ are the parameters more relevant in the irregular surface of tumor and $\alpha$ is the most important parameter for total area in the tumor growth.

\section{Discussion}\label{conlcusion}

In this paper we have presented a differential system for modelling the GBM growth for which we capture two properties according this kind of brain tumor: the ring width and the regularity of the tumor surface. 
\\

In order to detect these phenomena, we have made a numerical study with respect to the parameters of the model. After the simulations and the results obtained, we have proved that the parameters more relevant according to the tumor growth are $\kappa_1$ and $\alpha$. 
\\

For the tumor ring, where the vasculature is uniformly distributed, the results show that $\alpha$ is the most relevant parameter as we can observe in Figs $\ref{Ring_dif_kappa1}$-$\ref{Ring_dif_beta}$. In the case of surface regularity, where the vasculature is non-uniformly distributed, the parameter which produce more irregularity in the tumor surface is $\kappa_1$, see Figs $\ref{crec_irreg_kappa1}$-$\ref{crec_irreg_beta2}$. 
\\

However, for the total area in the surface regularity, the parameter $\alpha$ achieves the highest area for $\alpha=10$. Furthermore, in tumor growth section, despite the areas for $\kappa_1=100$ and $\alpha=10$ seem similar, the critical value from which the tumor is considered, defined in $\left(\ref{Tmin}\right)$, occupies more space for $\alpha=10$ than for $\kappa_1$ as we can see in Figs $\ref{alpha_10_e}$ and $\ref{kappa1_100_e}$. Hence, we can conclude that not only is $\alpha$ the main parameter for the tumor ring, but also it can increase or decrease the amount of total area with higher influence than $\kappa_1$.
\\

Finally, we have reduced our study from $9$ initial parameters to $2$ essentials parameters which determine the both main issues of GBM; the different tumor rings and the regular or irregular tumor surface. We have showed that $\alpha$ is the most relevant parameter related to the density and area of tumor independently the distribution of the vasculature. 
\addcontentsline{toc}{section}{References}
%\bibliographystyle{ws-m3as}
%\bibliography{reference.bib}

\begin{thebibliography}{9}
	\bibitem{Alfonso_2017}
	J. C. L. Alfonso et al., The biology and mathematical modelling of glioma invasion: a review. J. R. Soc. Interface. $\mathbf{14}$ (2017) 20170490. \url{https://doi.org/10.1098/rsif.2017.0490}.
	
	\bibitem{Baldock_2013}
	A. Baldock et al., From patient-specific mathematical neuro-oncology to precision medicine. Front. Oncol. $\mathbf{3}$ (2013) 62. \url{ https://doi.org/10.3389/fonc.2013.00062}.
	
	\bibitem{Davis_2016}
	M. E. Davis, Glioblastoma: Overview of disease and treatment, Clin. J. Oncol. Nurs., $\mathbf{20}$
	(2016), S2-S8. \url{10.1188/16.CJON.S1.2-8}.
	
	\bibitem{Romero2_2020}
	A. Fernández-Romero, F. Guillén-González and A. Suárez, Theoretical and numerical analysis for a hybrid tumor model with diffusion depending on vasculature. J. Math. Anal. Appl. $\mathbf{503}$ (2021) 29. \url{https://doi.org/10.1016/j.jmaa.2021.125325}.
	
	\bibitem{Klank_2018}
	R. L. Klank, S.S. Rosenfeld and D.J. Odde, A Brownian dynamics tumor progression simulator with application to glioblastoma. Converg. Sci. Phys. Oncol. $\mathbf{4}$ (2018) 015001. \url{10.1088/2057-1739/aa9e6e}.
	
	\bibitem{Alicia_2015}
	A. Martínez-González et al,. Combined therapies of antithrombotics and antioxidants delay in silico brain tumour progression. Math. Med. Biol. $\mathbf{32}$ (2015) 239-262. \url{https://doi.org/10.1093/imammb/dqu002}.
	
	\bibitem{Alicia_2012}
	A. Martínez-González, G. F. Calvo, L. A. Pérez-Romasanta and V.M. Pérez-García, Hypoxic cell waves around necrotic cores in glioblastoma: a mathematical model and its therapeutical implications. Bull. Math. Biol. $\mathbf{74}$ (2012) 2875-2896. \url{https://doi.org/10.1007/s11538-012-9786-1}.
	
	\bibitem{Ostrom_2014}
	Q. T. Ostrom et al., \uppercase{C}\uppercase{B}\uppercase{T}\uppercase{R}\uppercase{U}\uppercase{S} statistical report: primary brain and central nervous system tumors diagnosed in the united states in 2007-2011. Neuro-Oncol. $\mathbf{16}$ (2014) iv1-iv63. \url{https://doi.org/10.1093/neuonc/nou223}.
	
	\bibitem{Victor_2020}
	J. Pérez-Beteta, J. Belmonte-Beitia and V. M. Pérez-García, Tumor width on T1-weighted MRI images of glioblastoma as a prognostic biomarker: a mathematical model, Math. Model. Nat. Phenom. $\mathbf{15}$ (2020) 10. \url{https://doi.org/10.1051/mmnp/2019022}.
	
	\bibitem{Julian_2016}
	J. Pérez-Beteta et al., Glioblastoma: does the pretreatment geometry matter? A postcontrast T1 MRI-based study. Eur. Radiol. $\mathbf{27}$ (2017) 163-169. \url{https://doi.org/10.1007/s00330-016-4453-9}.
	
	\bibitem{Victor_2018}
	J. Pérez-Beteta et al., Tumor surface regularity at MR imaging predicts survival and response to surgery in patients with glioblastoma, Radiology. $\mathbf{288}$ (2018) 218-225. \url{https://doi.org/10.1148/radiol.2018171051}.
	
	\bibitem{Protopapa_2018}
	M. Protopapa et al., Clinical implications of in silico mathematical modeling for glioblastoma: a critical review, J. Neurooncol. $\mathbf{136}$ (2018) 1-11. \url{https://doi.org/10.1007/s11060-017-2650-2}.
	
	\bibitem{Rockne_2009}
	R. Rockne, E. Alvord, J. Rockhill, and K. Swanson, A mathematical model for brain tumor response to radiation therapy, J. Math. Biol., $\mathbf{58}$ (2009), 561-578. \url{https://doi.org/10.1007/s00285-008-0219-6}.
	
	\bibitem{Swanson_2008}
	K. Swanson, R. Rostomily, and E. Alvord, Jr., A mathematical modelling tool for predicting survival of individual patients following resection of glioblastoma: A proof of principle, British J. Cancer, $\mathbf{98}$ (2008), 113-119. \url{https://doi.org/10.1038/sj.bjc.6604125}.
	
	\bibitem{Swanson_2000}
	K. R. Swanson, E. C. Alvord, Jr., and J. Murray, A quantitative model for differential motility of gliomas in grey and white matter, Cell Prolif, $\mathbf{33}$ (2000), 317-329. \url{10.1046/j.1365-2184.2000.00177.x}.
	
	\bibitem{Tracqui_1995}
	P. Tracqui et al., A mathematical model of glioma growth: The effect of chemotherapy on spatio-temporal growth, Cell Prolif, $\mathbf{28}$ (1995), 17-31. \url{https://doi.org/10.1111/j.1365-2184.1995.tb00036.x}.
	
	
	
	
\end{thebibliography}

\end{document}